\documentclass[a4paper,12pt]{article} 

\usepackage{color}

\usepackage{amsfonts, amsmath, amsthm, amssymb}
\usepackage[T1]{fontenc}
\usepackage[cp1250]{inputenc}
\usepackage{xcolor}
\usepackage{graphicx}
\usepackage{amssymb}
\usepackage{amsmath}
\usepackage{mathptmx}
\usepackage{helvet}
\usepackage{courier}
\usepackage{txfonts}
\usepackage{tikz} 
\usetikzlibrary{arrows}
\usepackage{type1cm}

\usepackage{verbatim}

\usepackage{graphicx}
\usepackage{epsfig,amscd,amssymb,amsxtra,amsmath,amsthm}
\usepackage{type1cm}
\usepackage[T1]{fontenc}
\usepackage{graphics}
\usepackage[mathscr]{eucal}
\usepackage[all]{xy}
\usepackage{amsmath,amscd}

\newtheorem{theorem}{Theorem}[section]

\newtheorem{definition}[theorem]{Definition}

\newtheorem{example}[theorem]{Example}

\newtheorem{corollary}[theorem]{Corollary}
\newtheorem{problem}[theorem]{Problem}
\newtheorem{observation}[theorem]{Observation}

\newcommand{\diam} {\mathop{\rm diam}\nolimits}

\begin{document}

\def\joinrel{\mkern-3mu}
\newcommand{\varproj}{\displaystyle \lim_{\multimapinv\joinrel-\joinrel-}}

\title{Specification in Mahavier Systems via Closed Relations}
\author{Iztok Bani\v c, Goran Erceg, Ivan Jeli\' c, Judy Kennedy}
\date{}

\maketitle

\begin{abstract}
We study two fundamental properties of topological dynamical systems, the specification property and the initial specification property, and explore their generalizations to the broader setting of CR-dynamical systems, where the dynamics are governed by closed relations rather than continuous functions. While these two properties are equivalent for many classical systems, we demonstrate that their generalizations to CR-dynamical systems often lead to distinct behaviors. Applying them to Mahavier dynamical systems, we introduce new specification-type properties. These generalized notions extend the classical theory and reveal rich structural differences in dynamical behavior. Moreover, each of the new properties reduces to the standard specification property when restricted to continuous functions. \end{abstract}
\-
\\
\noindent
{\it Keywords:} Dynamical systems; Mahavier dynamical systems; closed relations, specification property\\
\noindent
{\it 2020 Mathematics Subject Classification:} 37B02, 37B45, 54C60, 54F15, 54F17

\section{Introduction}

In many cases, when constructing models for empirical data, continuous functions prove insufficient, and the data are more accurately represented by closed relations. A notable example of this phenomenon appears in macroeconomics through the Christiano-Harrison model (see \cite{model} for more details).

This illustrates the need to study closed relations on compact metric spaces; continuous functions simply cannot account for every situation. Since every continuous function is, in fact, a special case of a closed relation on a compact metric space, it is natural to generalize the classical notion of a dynamical system $(X, f)$ to what we call a CR-dynamical system, that is, a dynamical system based on a closed relation.

CR-dynamical systems thus provide a broader framework than traditional ones and are particularly useful in contexts where modeling with continuous functions falls short. Moreover, when we construct a closed relation $F$ on a space $X$, we are not just defining a CR-dynamical system $(X, F)$; we are simultaneously giving rise to two associated dynamical systems: the forward system $(X_F^+, \sigma_F^+)$ and the two-sided system $(X_F, \sigma_F)$, each offering further insights into the dynamics generated by $F$. They are called Mahavier dynamical systems.
		 
In this paper, we investigate two key properties of topological dynamical systems: the specification property and the initial specification property. For many classical systems, these two notions are equivalent, and as a result, the initial specification property is often used as a tool to study the specification property. However, there are known examples where the two properties diverge.

Despite their apparent similarity, these properties give rise to fundamentally different behaviors when extended to the broader setting of CR-dynamical systems and subsequently applied to Mahavier dynamical systems. Our aim is to generalize both the specification property and the initial specification property from topological dynamical systems to CR-dynamical systems, and then to analyze how these generalized notions differ within the class of Mahavier dynamical systems.

These generalizations turn out to be quite meaningful: each gives rise to a distinct specification-like property on Mahavier dynamical systems. Moreover, each of the newly defined properties reduces to the classical (initial) specification property when the family of closed relations is restricted to continuous functions.

A closely related approach to specification in dynamics was taken by Raines and Tennant in \cite{raines}, where they developed a notion of the specification property for upper semi-continuous set-valued functions $F: X \to 2^X$ on compact metric spaces. Their definition, formulated in terms of tracing orbit segments by periodic orbits, ensures strong dynamical applications such as topological mixing and positive topological entropy. They also proved that their variant of the specification property for $F$ implies the specification property for the associated inverse limit system $(\varprojlim F, \sigma)$, where $\sigma$ is the standard shift map. 

In contrast, our work is set in the more general framework of closed relations rather than set-valued functions, which allows us to capture a strictly larger class of dynamical systems. While every upper semi-continuous function defines a closed relation, the converse does not hold, and many natural examples in applications, including in economic modeling, arise in this broader setting. Furthermore, we develop two distinct generalizations of the classical specification property - one based on set-wise tracing of orbit segments ($\mathcal{SP}$), and the other defined via the Hausdorff metric on the hyperspace of closed sets ($\mathcal{HSP}$). These variants coincide in classical settings but may differ in CR-dynamical systems. However, the variant of the specification property from \cite{raines} does not generalize the classical specification property. Our definitions also naturally extend to Mahavier dynamical systems, enabling us to identify and analyze multiple types of specification-like properties.

We organize the paper as follows. In Section \ref{s1}, we recall the classical definition of the specification property in topological dynamical systems and establish some of its fundamental properties. We also introduce the initial specification property, which serves as a useful tool in the study of the specification property. Although these two concepts are typically equivalent in the classical setting, we will later demonstrate that their CR-dynamical generalizations, when applied to Mahavier systems, often differ. In Section \ref{s5}, Mahavier dynamical systems are defined and some general properties on such systems are presented. 
In Section \ref{s2}, we generalize the specification property to CR-dynamical systems and use this framework to define new variants of the specification property for Mahavier dynamical systems. In Section \ref{s3}, we carry out a parallel generalization of the initial specification property, again obtaining another corresponding variant for Mahavier dynamical systems. Section \ref{s4} is devoted to the study of the relationships and differences between these newly defined variants. 
\section{(Initial) specification property}\label{s1}
Since the specification property involves a relatively intricate definition, we begin this section by presenting it in detail. Before doing so, we introduce some fundamental concepts related to topological dynamical systems that are necessary for understanding the definition.

The specification property is formally defined in Definition \ref{def1}, but to fully grasp it, several preliminary definitions must first be given.
\begin{definition}
	Let $X$ be a compact metric space and let $f:X\rightarrow X$ be a continuous function. 
	We say that $(X,f)$ is \emph{a (topological) dynamical system}. 
\end{definition}

\begin{definition}
Let $(X,f)$ be a dynamical system, let $x\in X$ and let $k,\ell$ be non-negative integers. If $k\leq \ell$, then we say that 
$$
f^{[k,\ell]}(x)=\Big(f^k(x),f^{k+1}(x),f^{k+2}(x),\ldots,f^{\ell}(x)\Big)
$$
is the \emph{$[k,\ell]$-orbit segment of the point $x$}.  
\end{definition}

\begin{definition}
Let $(X,f)$ be a dynamical system, let $n$ be a positive integer and for each $j\in \{1,2,3,\ldots,n\}$, let 
\begin{enumerate}
	\item $k_j$ and $\ell_j$ be non-negative integers such that $k_j\leq \ell_j$, and
	\item $x_j\in X$.
\end{enumerate}
We say that the $n-tuple$
$$
\Bigg(f^{[k_1,\ell_1]}(x_1),f^{[k_2,\ell_2]}(x_2),f^{[k_3,\ell_3]}(x_3),\ldots ,f^{[k_n,\ell_n]}(x_n)\Bigg)
$$
is \emph{an $n$-specification} or just \emph{a specification in $(X,f)$}.
\end{definition}
\begin{definition}
Let $(X,f)$ be a dynamical system, let $N$ be a positive integer and let 
$$
\mathcal S=\Bigg(f^{[k_1,\ell_1]}(x_1),f^{[k_2,\ell_2]}(x_2),f^{[k_3,\ell_3]}(x_3),\ldots ,f^{[k_n,\ell_n]}(x_n)\Bigg)
$$
be a specification in $(X,f)$. We say that $\mathcal S$ is \emph{an $N$-spaced specification}, if for each $j\in\{1,2,3,\ldots,n-1\}$,
$$
k_{j+1}-\ell_j\geq N.
$$
\end{definition}
\begin{definition}
Let $(X,f)$ be a dynamical system, let $N$ be a positive integer, let $\varepsilon>0$, let $y\in X$ and let 
$$
\mathcal S=\Bigg(f^{[k_1,\ell_1]}(x_1),f^{[k_2,\ell_2]}(x_2),f^{[k_3,\ell_3]}(x_3),\ldots ,f^{[k_n,\ell_n]}(x_n)\Bigg)
$$
be an $N$-spaced specification in $(X,f)$. We say that \emph{$\mathcal S$ is $\varepsilon$-traced in $(X,f)$ by $y$} if for each $i\in \{1,2,3,\ldots,n\}$ and for each $j\in \{k_i,k_i+1,k_i+2,\ldots,\ell_i\}$,
$$
d(f^j(y),f^j(x_i))\leq \varepsilon.
$$
\end{definition}
Finally, in Definition \ref{def1}, the specification property is defined.
\begin{definition}\label{def1}
Let $(X,f)$ be a dynamical system. We say that \emph{$(X,f)$ has the specification property} if for each $\varepsilon >0$, there is a positive integer $N$ such that for any $N$-spaced specification $\mathcal S$ in $(X,f)$, there is $y\in X$ such that $\mathcal S$ is $\varepsilon$-traced in $(X,f)$ by $y$.
\end{definition}
In Definition \ref{def2}, the initial specification property is defined. The following definitions are needed to properly define it.
\begin{definition}
Let $(X,f)$ be a dynamical system and let 
$$
\mathcal S=\Bigg(f^{[k_1,\ell_1]}(x_1),f^{[k_2,\ell_2]}(x_2),f^{[k_3,\ell_3]}(x_3),\ldots ,f^{[k_n,\ell_n]}(x_n)\Bigg)
$$
be a specification in $(X,f)$. We say that $\mathcal S$ is \emph{an $(\ell_1,\ell_2,\ell_3,\ldots,\ell_n)$-specification} if for each $i\in\{1,2,3,\ldots,n\}$, $k_i=0$. 
\end{definition}
\begin{definition}
Let $(X,f)$ be a dynamical system and let $\mathcal S$ be a specification in $(X,f)$. We say that $\mathcal S$ is \emph{an initial specification}  if there are non-negative integers $\ell_1,\ell_2,\ell_3,\ldots,\ell_n$ such that $\mathcal S$ is an $(\ell_1,\ell_2,\ell_3,\ldots,\ell_n)$-specification.
\end{definition}
\begin{definition}
Let $(X,f)$ be a dynamical system, let $m_1,m_2,m_3,\ldots,m_{n-1}$ be positive integers, let $\varepsilon>0$, let $y\in X$, and let 
$$
\mathcal S=\Bigg(f^{[0,\ell_1]}(x_1),f^{[0,\ell_2]}(x_2),f^{[0,\ell_3]}(x_3),\ldots ,f^{[0,\ell_n]}(x_n)\Bigg)
$$
be an initial specification in $(X,f)$. We say that $\mathcal S$ is \emph{$(\varepsilon,m_1,m_2,m_3,\ldots,m_{n-1})$-traced in $(X,f)$ by $y$} if for each $i\in \{1,2,3,\ldots,n\}$ and for each $j\in \{0,1,2,\ldots,\ell_i\}$,
$$
d\Big(f^{\ell_1+m_1+\ell_2+m_2+\ell_3+m_3+\ldots+\ell_{i-1}+m_{i-1}+j}(y),f^j(x_i)\Big)\leq \varepsilon.
$$
\end{definition}
\begin{definition}\label{def2}
Let $(X,f)$ be a dynamical system. We say that $(X,f)$ has \emph{the initial specification property} if for each $\varepsilon >0$, there is a positive integer $N$ such that 
\begin{enumerate}
	\item for each positive integer $n$, 
	\item for all positive integers $m_1,m_2,m_3,\ldots,m_{n-1}$ such that for each $i\in \{1,2,3,\ldots,n\}$, $m_i\geq N$, and 
	\item for each initial $n$-specification $\mathcal S$,
\end{enumerate} 
 there is $y\in X$ such that $\mathcal S$ is $(\varepsilon,m_1,m_2,m_3,\ldots,m_{n-1})$-traced in $(X,f)$ by $y$. 
\end{definition}
In the following example, it is shown that, in general, the specification property and the initial specification property are not equivalent.
\begin{example}
	Let $X=[0,1]$ and let $f:X\rightarrow X$ be defined by
	$$
	f(x)=0
	$$
	for any $x\in X$. Note that $(X,f)$ has the specification property but it does not have the initial specification property.
\end{example}

The following theorem is mathematical folklore. It is well-known. For the sake of completeness, we state it and give its detailed proof. 

\begin{theorem}\label{lala}
	Let $X$ be a compact metric space and let $f:X\rightarrow X$ be a continuous surjection. The following statements are equivalent.
	\begin{enumerate}
		\item \label{1e} The dynamical system $(X,f)$ has the specification property.
		\item \label{2e} The dynamical system $(X,f)$ has the initial specification property.
	\end{enumerate} 
\end{theorem}
\begin{proof} First, we prove that \ref{1e} implies \ref{2e}.
	Suppose that $(X,f)$ has the specification property and let $\varepsilon>0$. Also, let $N$ be a positive integer such that for any $N$-spaced specification $\mathcal S$ in $(X,f)$, there is $y\in X$ such that $\mathcal S$ is $\varepsilon$-traced in $(X,f)$ by $y$. Such an $N$ does exist since $(X,f)$ has the specification property. We prove that for each positive integer $n$, for all positive integers $m_1,m_2,m_3,\ldots,m_{n-1}$ such that for each $i\in \{1,2,3,\ldots,n\}$, $m_i\geq N$, and for each initial $n$-specification $\mathcal S$, there is $y\in X$ such that $\mathcal S$ is $(\varepsilon,m_1,m_2,m_3,\ldots,m_{n-1})$-traced in $(X,f)$ by $y$. 

So, let $n$ be a positive integer, let $m_1,m_2,m_3,\ldots,m_{n-1}$ be positive integers such that for each $i\in \{1,2,3,\ldots,n\}$, $m_i\geq N$, and let 
$$
\mathcal S=\Bigg(f^{[0,\ell_1]}(x_1),f^{[0,\ell_2]}(x_2),f^{[0,\ell_3]}(x_3),\ldots ,f^{[0,\ell_n]}(x_n)\Bigg)
$$
 be an initial $n$-specification. Let $a_1=0$, let $b_1=\ell_1$, and for each $i\in \{2,3,4\ldots,n\}$, let 
\begin{enumerate}
	\item $x_i^1\in f^{-1}(x_i)$,
	\item $\sum_{k=1}^{i-1}(\ell_k+m_k)=a_i$,
	\item $\sum_{k=1}^{i-1}(\ell_k+m_k)+\ell_{i}=b_i$, and
	\item for each $j\in \{1,2,3,\ldots,a_i-1\}$, let $x_i^{j+1}\in f^{-1}(x_i^{j})$.
\end{enumerate} 
 Then 
 $$
\mathcal C=\Bigg(f^{[a_1,b_1]}(x_1),f^{[a_2,b_2]}(x_2^{a_2}),f^{[a_3,b_3]}(x_3^{a_3}),\ldots ,f^{[a_n,b_n]}(x_n^{a_n})\Bigg)
 $$
 is an $N$-spaced specification. Let $y\in X$ be such that $\mathcal C$ is $\varepsilon$-traced in $(X,f)$ by $y$. We claim that $\mathcal S$ is $(\varepsilon,m_1,m_2,m_3,\ldots,m_{n-1})$-traced in $(X,f)$ by $y$. 
Let $i\in \{1,2,3,\ldots,n\}$ and let $j\in \{0,1,2,\ldots,\ell_i\}$. Then 
\begin{align*}
	&d\Big(f^{\ell_1+m_1+\ell_2+m_2+\ell_3+m_3+\ldots+\ell_{i-1}+m_{i-1}+j}(y),f^j(x_i)\Big)=\\
	&d\Big(f^{\ell_1+m_1+\ell_2+m_2+\ell_3+m_3+\ldots+\ell_{i-1}+m_{i-1}+j}(y),f^{\ell_1+m_1+\ell_2+m_2+\ell_3+m_3+\ldots+\ell_{i-1}+m_{i-1}+j}(x_i^{a_i})\Big)\leq \varepsilon.
\end{align*}
This proves that \ref{1e} implies \ref{2e}. To prove that \ref{2e} implies \ref{1e}, suppose that $(X,f)$ has the initial specification property and let $\varepsilon > 0$. Also, let $N$ be a positive integer such that for each positive integer $n$, for all positive integers $m_1,m_2,m_3,\ldots,m_{n-1}$ such that for each $i\in \{1,2,3,\ldots,n\}$, $m_i\geq N$, and for each initial $n$-specification $\mathcal S$, there is $y\in X$ such that $\mathcal S$ is $(\varepsilon,m_1,m_2,m_3,\ldots,m_{n-1})$-traced in $(X,f)$ by $y$. Such $N$ does exist since  $(X,f)$ has the initial specification property. We prove that for any $N$-spaced specification $\mathcal S$ in $(X,f)$, there is $y\in X$ such that $\mathcal S$ is $\varepsilon$-traced in $(X,f)$ by $y$. Let 
$$
\mathcal S=\Bigg(f^{[k_1,\ell_1]}(x_1),f^{[k_2,\ell_2]}(x_2),f^{[k_3,\ell_3]}(x_3),\ldots ,f^{[k_n,\ell_n]}(x_n)\Bigg)
$$
be an $N$-spaced specification in $(X,f)$. Then 
$$
\mathcal C=\Bigg(f^{[0,\ell_1-k_1]}(f^{k_1}(x_1)),f^{[0,\ell_2-k_2]}(f^{k_2}(x_2)),f^{[0,\ell_3-k_3]}(f^{k_3}(x_3)),\ldots ,f^{[0,\ell_n-k_n]}(f^{k_n}(x_n))\Bigg)
$$
is an initial $n$-specification. For each $i\in \{1,2,3,\ldots, n-1\}$, let $m_i=k_{i+1}-\ell_i$. Note that for each $i\in \{1,2,3,\ldots, n-1\}$, $m_i\geq N$ holds. Let $z\in X$ be such that $\mathcal C$ is $(\varepsilon,m_1,m_2,m_3,\ldots,m_{n-1})$-traced in $(X,f)$ by $z$. Also, let $y\in X$ be such that $z=f^{k_1}(y)$. We claim that $\mathcal S$ is $\varepsilon$-traced in $(X,f)$ by $y$.
Let $i\in \{1,2,3,\ldots,n\}$ and let $j\in \{k_i,k_i+1,k_i+2,\ldots,\ell_i\}$. Then 
\begin{align*}
&d(f^j(y),f^j(x_i))=\\
&d(f^{(\ell_1-k_1)+m_1+(\ell_2-k_2)+m_2+(\ell_3-k_3)+m_3+\ldots+(\ell_{i-1}-k_{i-1})+m_{i-1}+(j-k_i)}(z),f^{j-k_i}(f^{k_i}(x_i)))\leq \varepsilon.
\end{align*}
This completes the proof. 
\end{proof}
Next, we prove in Theorem \ref{lalalala} that both the specification property and  the initial specification property are dynamical properties. 
\begin{definition}
	Let $(X,f)$ and $(Y,g)$ be dynamical systems. If there is a homeomorphism $\varphi:X\rightarrow Y$ such that 
$$
\varphi \circ f=g\circ \varphi,
$$
then we say that \emph{$(X,f)$ and $(Y,g)$ are topological conjugates}.
\end{definition}
\begin{observation}
	Let $(X,f)$ and $(Y,g)$ be dynamical systems and let $\varphi:X\rightarrow Y$ be a homeomorphism such that $\varphi \circ f=g\circ \varphi$. Then for each positive integer $n$,
	$$
	g^n=\varphi\circ f^n\circ \varphi^{-1}.
	$$
\end{observation}
\begin{theorem}\label{lalalala}
	Let $(X,f)$  and $(Y,g)$ be dynamical systems and suppose that $(X,f)$ has the (initial) specification property. If $(X,f)$ and $(Y,g)$ are topological conjugates, then also $(Y,g)$ has the (initial) specification property.
\end{theorem}
\begin{proof}
First, let $d_X$ be the metric on $X$ and $d_Y$ the metric on $Y$. Suppose that $(X,f)$ and $(Y,g)$ are topological conjugates and let $\varphi:X\rightarrow Y$ be a homeomorphism such that $\varphi\circ f=g\circ \varphi$. First, suppose that $(X,f)$ has the specification property. We prove  that $(Y,g)$ has the specification property by showing that for each $\varepsilon >0$, there is a positive integer $N$ such that for any $N$-spaced specification $\mathcal S$ in $(Y,g)$, there is $y\in Y$ such that $\mathcal S$ is $\varepsilon$-traced in $(Y,g)$ by $y$. Let $\varepsilon >0$ and let $\delta>0$ be such that for all $x_1,x_2\in X$,
$$
d_X(x_1,x_2)<\delta ~~~ \Longrightarrow ~~~ d_Y(\varphi(x_1),\varphi(x_2))<\varepsilon.
$$ 
Since $(X,f)$ has the specification property, there is a positive integer $N$ such that for any $N$-spaced specification $\mathcal S$ in $(X,f)$, there is $x\in X$ such that $\mathcal S$ is $\frac{\delta}{2}$-traced in $(X,f)$ by $x$. Choose and fix such a positive integer $N$ and let 
$$
\mathcal S=\Bigg(g^{[k_1,\ell_1]}(y_1),g^{[k_2,\ell_2]}(y_2),g^{[k_3,\ell_3]}(y_3),\ldots ,g^{[k_n,\ell_n]}(y_n)\Bigg)
$$
 be an $N$-spaced specification in $(Y,g)$.  Then 
 $$
\mathcal C=\Bigg(f^{[k_1,\ell_1]}(\varphi^{-1}(y_1)),f^{[k_2,\ell_2]}(\varphi^{-1}(y_2)),f^{[k_3,\ell_3]}(\varphi^{-1}(y_3)),\ldots ,f^{[k_n,\ell_n]}(\varphi^{-1}(y_n))\Bigg)
$$
 is an $N$-spaced specification in $(X,f)$. Let $x$ be such that $\mathcal C$ is $\frac{\delta}{2}$-traced in $(X,f)$ by $x$ and let $y=\varphi(x)$. To prove that $\mathcal S$ is $\varepsilon$-traced in $(Y,g)$ by $y$, let $i\in \{1,2,3,\ldots,n\}$ and let $j\in \{k_i,k_i+1,k_i+2,\ldots,\ell_i\}$. We prove that $d_Y(g^j(y),g^j(y_i))\leq \varepsilon$. First, note that $d_X(f^j(x),f^j(\varphi^{-1}(y_i)))<\delta$. Therefore, $d_Y(\varphi(f^j(x)),\varphi(f^j(\varphi^{-1}(y_i))))<\varepsilon$. It follows from 
 $$
d_Y(g^j(y),g^j(y_i))=d_Y(\varphi(f^j(\varphi^{-1}(y))),\varphi(f^j(\varphi^{-1}(y_i))))=d_Y(\varphi(f^j(x)),\varphi(f^j(\varphi^{-1}(y_i))))
 $$ 
 that $d_Y(g^j(y),g^j(y_i))\leq \varepsilon$. The second part of the proof for the initial specification property is analogous to the first part. We leave it for the reader.
\end{proof}
\begin{definition}
	Let $\mathcal P$ be a property. We say that $\mathcal P$ is \emph{a dynamical property}, if for any two dynamical systems $(X,f)$ and $(Y,g)$ that are topological conjugates, the following holds:
	$$
	(X,f) \textup{ has property } \mathcal P ~~~ \Longrightarrow ~~~ (Y,g) \textup{ has property } \mathcal P.
	$$
\end{definition}
\begin{observation}
	Theorem \ref{lalalala} says that the specification property and the initial specification property are dynamical properties.
\end{observation}
For a dynamical property $\mathcal P$, the property $\mathcal P$ is often carried from the dynamical system $(X,f)$ to the dynamical system $\left(\varprojlim(X,f),\sigma\right)$, where $\sigma$ is the shift mapping on $\varprojlim(X,f)$. Theorem \ref{ivan} shows that this is also the case for the specification property and the initial specification property. First, the definition of an inverse limit is given. 
\begin{definition}
	Let $X$ be a compact metric space and let $f:X\rightarrow X$ be a continuous function. The \emph{inverse limit} generated by $(X,f)$ is the subspace
\begin{equation*}
 \varprojlim(X,f)=\Big\{(x_{1},x_{2},x_{3},\ldots ) \in \prod_{i=1}^{\infty} X \ | \ 
\text{ for each positive integer } i,x_{i}= f(x_{i+1})\Big\}
\end{equation*}
of the topological product $\prod_{i=1}^{\infty} X$.  { The function  $\sigma : \varprojlim(X,f) \rightarrow \varprojlim(X,f)$, 
 defined by 
$$
\sigma (x_1,x_2,x_3,x_4,\ldots )=(x_2,x_3,x_4,\ldots )
$$
for each $(x_1,x_2,x_3,\ldots )\in \varprojlim(X,f)$, 
is called \emph{   the shift map on $\varprojlim(X,f)$}.    }
\end{definition}
\begin{observation}\label{judo}
	Note that the shift map $\sigma$ on the inverse limit $\varprojlim(X,f)$ is a homeomorphism. Also, note that for each $(x_1,x_2,x_3,\ldots )\in \varprojlim(X,f)$,
	$$
	\sigma^{-1} (x_1,x_2,x_3,\ldots )=(f(x_1),x_1,x_2,x_3,\ldots ).
	$$
\end{observation}
\begin{theorem}\label{ivan}
	Let $\left(X,f\right)$ be a dynamical system, let $X_{\infty}=\varprojlim\left(X,f\right)$, and let $\sigma:X_{\infty}\rightarrow X_{\infty}$ be the shift map  on $X_{\infty}$. If $f$ is surjective, then the following statements are equivalent.
	\begin{enumerate}
		\item\label{g1} The dynamical system $\left(X,f\right)$ has the specification property 
		\item\label{g2} The dynamical system $\left(X_{\infty},\sigma\right)$ has the specification property.
		\item\label{g11} The dynamical system $\left(X,f\right)$ has the initial specification property 
		\item\label{g22} The dynamical system $\left(X_{\infty},\sigma\right)$ has the initial specification property.
	\end{enumerate} 
\end{theorem}	

\begin{proof}
Without loss of generality, we assume that $\diam(X)\leq 1$. To prove the implication from \ref{g1} to \ref{g2}, suppose that $\left(X,f\right)$ has the specification property and let $\varepsilon>0$. Let $j_0$ be a positive integer such that $\frac{1}{2^{j_0}}<\varepsilon$ and let $N_f$ be a positive integer such that for each $N_f$-spaced specification $\mathcal{S}_f$ in $\left(X,f\right)$ there exists $y\in X$ such that $\mathcal{S}_f$ is $\varepsilon$-traced in $(X,f)$ by $y$. 
	Put $N_{\sigma}=N_f+j_0-1$. Then $N_{\sigma}$ is a positive integer. Also, let $\mathbf{x_1}=\left(x^1_1,x^1_2,x^1_3,\ldots\right)$, $\mathbf{x_2}=\left(x^2_1,x^2_2,x^2_3,\ldots\right)$, $\mathbf{x_3}=\left(x^3_1,x^3_2,x^3_3,\ldots\right)$, $\ldots$, $\mathbf{x_n}=\left(x^n_1, x_2^n,x_3^n,\ldots\right)$ be any points in $X_{\infty}$ and let 
	$$
	\mathcal{S}_{\sigma}=\left(\sigma^{\left[k_1,\ell_1\right]}\left(\mathbf{x_1}\right),\sigma^{\left[k_2,\ell_2\right]}\left(\mathbf{x_2}\right),\sigma^{\left[k_3,\ell_3\right]}\left(\mathbf{x_3}\right),\ldots,\sigma^{\left[k_n,\ell_n\right]}\left(\mathbf{x_n}\right)\right),
	$$
	be an $N_{\sigma}$-spaced specification in $\left(X_{\infty},\sigma\right)$. Also, let $n_0=\ell_n+j_0+1$ and for each $i\in \{1,2,3,\ldots,n\}$, let 
	$$
	x_i=x^{n-i+1}_{n_0}.
	$$
	For each $i\in \{1,2,3,\ldots, n\}$, let  $k'_i=n_0-\ell_{n-i+1}-j_0$ and $\ell'_i=n_0-k_{n-i+1}-1$, and let 
	 $$
	\mathcal{S}_{f}=\left( f^{\left[k'_1,\ell'_1\right]}\left(x_1\right),f^{\left[k'_2,\ell'_2\right]}\left(x_2\right),f^{\left[k'_3,\ell'_3\right]}\left(x_3\right),\ldots,f^{\left[k'_n,\ell'_n\right]}\left(x_n\right)\right)
	$$ 
	be a specification in $\left(X,f\right)$. Since for arbitrary $i\in\left\lbrace1,2,3,\ldots,n-1\right\rbrace$,  it holds that
	$$
	k'_{i+1}-\ell'_i=n_0-\ell_{n-\left(i+1\right)+1}-j_0-n_0+k_{n-i+1}+1=
	$$
	$$
	=k_{n-i+1}-\ell_{n-i}-j_0+1\geq N_{\sigma}-j_0+1=N_f,
	$$
	it follows that $\mathcal{S}_{f}$ is an $N_f$-spaced specification in $\left(X,f\right)$.
	Hence, there exists $y\in X$ such that $\mathcal{S}_f$ is $\varepsilon$-traced in $(X,f)$ by $y$.
	Next, let 
	$$
	\mathbf{y}=\left(y_1,y_2,y_3,\ldots,y_{n_0-1},y_{n_0},y_{n_0+1},\ldots\right)\in X_{\infty}
	$$
	be such that $y_{n_0}=y$. 
	Such an element $\mathbf{y}$ in $X_{\infty}$ exists since $f$ is a surjection. We prove that $\mathcal{S}_{\sigma}$ is $\varepsilon$-traced in $(X_{\infty},\sigma)$ by $\mathbf{y}$. Let $i\in\left\lbrace 1,2,3,\ldots,n\right\rbrace$ and take an arbitrary $j\in\left\lbrace k_i,k_i+1,k_i+2,\ldots,\ell_i\right\rbrace$. Also, let 
\begin{align*}
	&M_1=\left\lbrace\frac{d\left(x^i_{j+m},y_{j+m}\right)}{2^m} \ \Big| \ m\in \{1,2,3,\ldots,\ell_i+j_0-j\}\right\rbrace\\
	&M_2=\left\lbrace\frac{d\left(x^i_{j+m},y_{j+m}\right)}{2^m} \ \Big| \ m\in \{\ell_i+j_0-j,\ell_i+j_0-j+1,\ell_i+j_0-j+2,\ldots\}\right\rbrace.
\end{align*}
	Then
	\begin{align*}
		&\textup{d}_{\sup}\left(\sigma^j\left(\mathbf{x_i}\right),\sigma^j\left(\mathbf{y}\right)\right)=\textup{d}_{\sup}\left(\left(x^i_{j+1},x^i_{j+2},x^i_{j+3},\ldots\right),\left(y_{j+1},y_{j+2},y_{j+3},\ldots\right)\right)=\\
		&\max\left\lbrace\frac{d\left(x^i_{j+m},y_{j+m}\right)}{2^m} \ \Big| \ m \ \textup{is a positive integer }\right\rbrace=\max\left\lbrace \max M_1,\max M_2\right\rbrace.
	\end{align*}
		Note that 
		$$
		j_0\leq \ell_i+j_0-j\leq \ell_i-k_i+j_0.
		$$
Since for each positive integer $m$, $d\left(x^i_{j+m},y_{j+m}\right)\leq 1$ and since $\ell_i+j_0-j\geq j_0$, it follows that 
	$$
	\max M_2\leq\frac{1}{2^{j_0+1}}<\varepsilon.
	$$
	Furthermore, for each $m\in \{1,2,3,\ldots,\ell_i+j_0-j\}$, it holds that
	\begin{align*}
		&d\left(x^i_{j+m},y_{j+m}\right)=d\left(f^{n_0-j-m}\left(x^i_{n_0}\right),f^{n_0-j-m}\left(y_{n_0}\right)\right)=\\
		&d\left(f^{n_0-j-m}\left(x_{n-i+1}\right),f^{n_0-j-m}\left(y\right)\right)<\varepsilon
	\end{align*}
		since
	$$
	n_0-j-m\geq n_0-j-\ell_i-j_0+j=n_0-\ell_i-j_0=k'_{n-i+1}\quad\textrm{and}
	$$
	$$
	n_0-j-m\leq n_0-j-1\leq n_0-k_i-1=\ell'_{n-i+1}.
	$$
	Thus,
	$$
	\max M_1<\varepsilon
	$$
	and, finally, it follows that
	$$
	\textup{d}_{\sup}\left(\sigma^j\left(\mathbf{x_i}\right),\sigma^j\left(\mathbf{y}\right)\right)<\varepsilon.
	$$
	It follows that $\left(X_{\infty},\sigma\right)$ has the specification property. This proves the implication from \ref{g1} to \ref{g2}. 
	
		To prove the implication from \ref{g2} to \ref{g1}, suppose that $\left(X_{\infty},\sigma\right)$ has the specification property and let $\varepsilon>0$. Then there exists $N_{\sigma}\in\mathbb{N}$ such that for each $N_{\sigma}$-spaced specification $\mathcal{S}_{\sigma}$ in $\left(X_{\infty},\sigma\right)$ there exists $\mathbf{y}\in X_{\infty}$ such that $\mathcal{S}_{\sigma}$ is $\frac{\varepsilon}{2}$-traced by $\mathbf{y}$. Choose and fix such a positive integer $N_{\sigma}$, let $N_f=N_{\sigma}$, and let
	$$
	\mathcal{S}_{f}=\left( f^{\left[k_1,\ell_1\right]}\left(x_1\right),f^{\left[k_2,\ell_2\right]}\left(x_3\right),f^{\left[k_3,\ell_3\right]}\left(x_2\right),\ldots, f^{\left[k_n,\ell_n\right]}\left(x_n\right)\right)
	$$
	be an arbitrary $N_f$-spaced specification in $\left(X,f\right)$. For each $i\in \{1,2,3,\ldots, n\}$, we choose and fix 
	$$
	\mathbf{x_i}=\left(x_1^i,x_2^i,x_3^i,\ldots,x^i_{\ell_n+1},x^i_{\ell_n+2},x^i_{l_n+3},\ldots\right) \in X_{\infty}
		$$
		such that $x^i_{\ell_n+2}=x_{n-i+1}$. Such elements $\mathbf{x_i}$ in $X_{\infty}$ exist since $f$ is a surjection. Now, let
	$$
	\mathcal{S}_{\sigma}=\left(\sigma^{\left[k'_1,\ell'_1\right]}\left(\mathbf{x_1}\right),\sigma^{\left[k'_2,\ell'_2\right]}\left(\mathbf{x_2}\right),\sigma^{\left[k'_3,\ell'_3\right]}\left(\mathbf{x_3}\right),\ldots,\sigma^{\left[k'_n,\ell'_n\right]}\left(\mathbf{x_n}\right)\right),
	$$
	where for each $i\in\{1,2,3,\ldots,n\}$,
	$$
	k'_i=\ell_n-\ell_{n-i+1}+1  ~~~ \textup{and} ~~~ \ell'_i=\ell_n-k_{n-i+1}+1.
	$$
	Since for an arbitrary $i\in\left\lbrace1,2,3,\ldots,n-1\right\rbrace$, 
	$$
	k'_{i+1}-\ell'_i=\ell_n-\ell_{n-\left(i+1\right)+1}+1-\ell_n+k_{n-i+1}-1=k_{n-i+1}-\ell_{n-i}\geq N_f=N_{\sigma},
	$$
	it follows that $\mathcal{S}_{\sigma}$ is $N_{\sigma}$-spaced specification in $\left(X,f\right)$. 	Hence, there exists $\mathbf{y}\in X_{\infty}$ such that $\mathcal{S}_{\sigma}$ is $\frac{\varepsilon}{2}$-traced in $(X_{\infty},\sigma)$ by $\mathbf{y}$. Choose and fix such a $\mathbf y=(y_1,y_2,y_3,\ldots)$
	and let $y=y_{\ell_n+2}$. We prove that $\mathcal{S}_f$ is $\varepsilon$-traced in $(X,f)$ by $y$. To do so, let $i\in\left\lbrace 1,2,3,\ldots,n\right\rbrace$. Then for an arbitrary $j\in\left\lbrace k'_{n-i+1},k'_{n-i+1}+1,k'_{n-i+1}+2,\ldots,\ell'_{n-i+1}\right\rbrace$, it holds that
	\begin{align*}
		\textup{d}_{\sup}\left(\sigma^j\left(\mathbf{x_{n-i+1}}\right),\sigma^j\left(\mathbf{y}\right)\right)=&\textup{d}_{\sup}\left(\left(x^{n-i+1}_{j+1},x^{n-i+1}_{j+2},x^{n-i+1}_{j+3},\ldots\right),\left(y_{j+1},y_{j+2},y_{j+3},\ldots\right)\right)	=\\
		&\max\left\lbrace\frac{d\left(x^{n-i+1}_{j+m},y_{j+m}\right)}{2^m} \ \Big| \  \ m \textup{ is a positive integer}\right\rbrace<\frac{\varepsilon}{2}.
	\end{align*}
	In particular, for $m=1$ and for every $j\in\left\lbrace k'_{n-i+1},k'_{n-i+1}+1,k'_{n-i+1}+2,\ldots,\ell'_{n-i+1}\right\rbrace$, it holds that
	$$
	\frac{d\left(x^{n-i+1}_{j+1},y_{j+1}\right)}{2}<\frac{\varepsilon}{2} \implies d\left(x^{n-i+1}_{j+1},y_{j+1}\right)<\varepsilon.
	$$
	Hence, for every $j\in\left\lbrace k_i,k_i+1,k_i+2,\ldots,\ell_i\right\rbrace$,
	$$
	d\left(f^j\left(x_i\right),f^j\left(y\right)\right)=d\left(f^j\left(x_{\ell_n+2}^{n-i+1}\right),f^j\left(y_{\ell_n+2}\right)\right)=d\left(x^{n-i+1}_{\ell_n+2-j},y_{\ell_n+2-j}\right)<\varepsilon
	$$
	since
	$$
	\ell_n+2-j\geq \ell_n+2-\ell_i=\ell_n-\ell_i+1+1=k'_{n-i+1}+1
	$$
	and
	$$
	\ell_n+2-j\leq \ell_n+2-k_i=\ell_n-k_i+1+1=\ell'_{n-i+1}+1.
	$$
	This proves that $\left(X,f\right)$ has the specification property. Note that wee have just proved that \ref{g1} is equivalent to \ref{g2}. Since $f$ is surjective, it follows from Theorem \ref{lala} that \ref{g11} is equivalent to \ref{g1}. Note that it follows from Observation \ref{judo} that $\sigma$ is surjective. Therefore, it follows from Theorem \ref{lala} that \ref{g22} is equivalent to \ref{g2}.
	\end{proof}
	\begin{observation}
		Let $\left(X,f\right)$ be a dynamical system, let $X_{\infty}=\varprojlim\left(X,f\right)$, and let $\sigma:X_{\infty}\rightarrow X_{\infty}$ be the shift map on $X_{\infty}$. Note that $\sigma$ is a surjection, even in the case where $f$ is not. Therefore, the following statements are equivalent.
	\begin{enumerate}
		\item\label{g222} The dynamical system $\left(X_{\infty},\sigma\right)$ has the specification property.
		\item\label{g221} The dynamical system $\left(X_{\infty},\sigma\right)$ has the initial specification property.
	\end{enumerate} 
	\end{observation}
	For a dynamical property $\mathcal P$, the property $\mathcal P$ is often carried from the dynamical system $(X,f)$, where $f$ is a homeomorphism on $X$, to the dynamical system $(X,f^{-1})$. Theorem \ref{goran} shows that this is also the case for the specification property and the initial specification property.
\begin{theorem}\label{goran}
	Let $(X,f)$ be a dynamical system. If $f$ is a homeomorphism, then the following statements are equivalent.
	\begin{enumerate}
		\item\label{juha1} The dynamical system $(X,f)$ has the specification property.
		\item\label{juha2} The dynamical system $\left(X,f^{-1}\right)$ has the specification property.
		\item\label{juha3} The dynamical system $(X,f)$ has the initial specification property.
		\item\label{juha4} The dynamical system $\left(X,f^{-1}\right)$ has the initial specification property.
	\end{enumerate}
\end{theorem}
\begin{proof}
Suppose that $(X,f)$ has the specification property and let $\varepsilon>0$. Let $N$ be a positive integer such that for any $N$-spaced specification $\mathcal S_f$ in $(X,f)$, there is $y\in X$ such that $\mathcal S_f$ is $\varepsilon$-traced in $(X,f)$ by $y$. We prove that for the same $N$, for any $N$-spaced specification $\mathcal S_{f^{-1}}$ in $(X,f^{-1})$ there is $y'\in X$ such that $\mathcal S_{f^{-1}}$ is $\varepsilon$-traced in $(X,f^{-1})$ by $y$.
	So let $n$ be a positive integer and let $$\mathcal S_{f^{-1}}=\left((f^{-1})^{[k_1,\ell_1]}(x_1),(f^{-1})^{[k_2,\ell_2]}(x_2),(f^{-1})^{[k_3,\ell_3]}(x_3),\ldots ,(f^{-1})^{[k_n,\ell_n]}(x_n)\right)$$ be any $N$-spaced specification in $(X,f^{-1})$. 
	We show that $\mathcal S_{f^{-1}}$ is an $N$-spaced specification in $(X,f)$ as follows. Denote 
	$$
	x'_1=(f^{-1})^{\ell_n}(x_n).
	$$ 
	Now we have 
	$$
	f(x'_1)=(f^{-1})^{\ell_n-1}(x_n), \ldots, f^{\ell_n-k_n}(x'_1)=(f^{-1})^{k_n}(x_n).
	$$ 
	Let $x'_2$ be a point in $X$ such that 
	$$
	f^{\ell_n-k_n+k_n-\ell_{n-1}}(x'_2)=f^{\ell_n-\ell_{n-1}}(x'_2)=(f^{-1})^{\ell_{n-1}}(x_{n-1}).
	$$
	Such a point $x'_2$ exists since $f$ is a homeomorphism.
	We have 
	$$
	f^{\ell_n-\ell_{n-1}+1}(x'_2)=(f^{-1})^{\ell_{n-1}-1}(x_{n-1}), \ldots, f^{\ell_n-k_{n-1}}(x'_2)= (f^{-1})^{k_{n-1}}(x_{n-1}).
	$$
	We proceed inductively and finally, let $x'_n$ be a point in $X$ such that
	$$
	f^{\ell_n-k_2+k_2-\ell_1}(x'_n)=f^{\ell_n-\ell_1}(x'_n)=(f^{-1})^{\ell_1}(x_1).
	$$ 
	Such a point $x'_n$ exists since $f$ is a homeomorphism.
	We have 
	$$
	f^{\ell_n-\ell_1+1}(x'_n)=(f^{-1})^{\ell_1-1}(x_1),\ldots,f^{\ell_n-k_1}(x'_n)=(f^{-1})^{k_1}(x_1).
	$$
	Let 
	$$
	\mathcal S_f=\left( f^{[0,\ell_n-k_n]}(x'_1), f^{[\ell_n-\ell_{n-1}, \ell_n-k_{n-1}]}(x'_2), \ldots, f^{[\ell_n-\ell_1,\ell_n-k_1 ]}(x'_n)  \right). 
	$$
	Since we defined all the orbit segments in $\mathcal S_{f^{-1}}$ from the ones in $\mathcal S_f$ just going backwards, $\mathcal S_f$ is obviously an  $N$-spaced specification in $(X,f)$.  
	Let $y \in X$ be a point such that $\mathcal S_f$ is $\varepsilon$-traced in $(X,f)$ by $y$.
	Let $y'$ be a point in $X$ such that $(f^{-1})^{k_1}(y')=f^{l_n-k_1}(y).$  Such a point $y'$ exists since $f^{-1}$ is a homeomorphism.
	Note that, since $\mathcal S_f$ is $\varepsilon$-traced in $(X,f)$ by $y$, we have that $\mathcal S_{f^{-1}}$ is $\varepsilon$-traced in $(X,f^{-1})$ by $y'.$  
	The other implication now follows by putting $f=g^{-1}$ and $f^{-1}=g$.
	
	Note that we have just proved that \ref{juha1} is equivalent to \ref{juha2}. Since $f$ is surjective, it follows from Theorem \ref{lala} that \ref{juha1} is equivalent to \ref{juha3}. Note that $f^{-1}$ is also surjective. Therefore, it follows from Theorem \ref{lala} that \ref{juha2} is equivalent to \ref{juha4}.

\end{proof}
\section{CR-dynamical systems and Mahavier dynamical systems}\label{s5}
In this section, we first define the notion of a CR-dynamical system. Building on this, we introduce the concept of a Mahavier dynamical system and then examine those Mahavier systems that exhibit the (initial) specification property. Along the way, we also define a new class of dynamical properties that arise naturally from the CR-dynamical setting.
\begin{definition}
	Let $X$ be a {non-empty} compact metric space and let $F\subseteq X\times X$ be a relation on $X$. If ${F}$ is closed in $X\times X$, then we say that ${F}$ is  \emph{a closed relation on $X$}.  Also, if $F$ is a closed relation on $X$, then we call $(X,F)$ \emph{a CR-dynamical system}.
\end{definition} 
\begin{observation}
	Let $(X,F)$ be a CR-dynamical system. Note that if for each $x\in X$ there is exactly one $y\in X$ such that $(x,y)\in F$, then $F$ is a continuous function from $X$ to $X$ and  \begin{enumerate}
		\item we write $F:X\rightarrow X$ instead of ${F}\subseteq X\times X$, and 
		\item for all $x,y\in X$, we write $y=F(x)$ instead of $(x,y)\in F$.
	\end{enumerate}
	\end{observation}
	 The following theorem is a well-known result, see \cite{akin} for more information.
	 \begin{theorem}
Let $(X,f)$ and $(Y,g)$ be dynamical systems. The following statements are equivalent.
\begin{enumerate}
\item The dynamical systems $(X,f)$ and $(Y,g)$ are topological conjugates. 
\item There is a homeomorphism $\varphi:X\rightarrow Y$  such that for each $(x,y)\in X\times X$, the following holds:
$$
(x,y)\in f \Longleftrightarrow (\varphi(x),\varphi(y))\in g.
$$
\end{enumerate}
\end{theorem}
	The following generalizes the notion of topological conjugates from topological dynamical systems to CR-dynamical systems.
	\begin{definition}
		Let $(X,F)$ and $(Y,G)$ be CR-dynamical systems. We say that \emph{$(X,F)$ and $(Y,G)$ are topological conjugates} if there is a homeomorphism $\varphi:X\rightarrow Y$ such that for each $(x,y)\in X\times X$, the following holds
 $$
 (x,y)\in F  \Longleftrightarrow (\varphi(x), \varphi(y))\in G.
 $$
	\end{definition}
	
\begin{theorem}\label{dekanicaMaja}
	Let $(X,F)$ and $(Y,G)$ be CR-dynamical systems that are topologically conjugate and let $\varphi:X\rightarrow Y$ be a homeomorphism such that for each $(x,y)\in X\times X$, the following holds
 $$
 (x,y)\in F  \Longleftrightarrow (\varphi(x), \varphi(y))\in G.
 $$
 Then for each $x,y\in X$,
 $$
 y\in F(x) ~~~ \Longrightarrow ~~~ \varphi(y)\in G(\varphi(x)).  
 $$
\end{theorem}
\begin{proof}
	Let $x,y\in X$ be such that $y\in F(x)$. Note that  $\varphi(y)\in \varphi(F(x))$. To see that $\varphi(y)\in G(\varphi(x))$ we prove that $\varphi(F(x))\subseteq G(\varphi(x))$. Let $z\in \varphi(F(x))$ and let $w\in F(x)$ be such that $z=\varphi(w)$. It follows that $(x,w)\in F$ and, therefore, $(\varphi(x),\varphi(w))\in G$. Hence, $z\in G(\varphi(x))$.
\end{proof}
\begin{definition}
	A property $\mathcal P$ is \emph{a CR-dynamical property}, if for any two CR-dynamical systems $(X,F)$ and $(Y,G)$ that are topological conjugates, the following holds:
	$$
	(X,F) \textup{ has property } \mathcal P ~~~ \Longrightarrow ~~~ (Y,G) \textup{ has property } \mathcal P.
	$$
\end{definition}	
\begin{observation}
	Let $\mathcal P$ be a CR-dynamical property and let   $(X,F)$ and $(Y,G)$ be CR-dynamical systems that are topological conjugates. Note that the following are equivalent: 
	\begin{enumerate}
		\item $(X,F) \textup{ has property } \mathcal P ~~~ \Longrightarrow ~~~ (Y,G) \textup{ has property } \mathcal P$.
		\item $(X,F) \textup{ has property } \mathcal P ~~~ \Longleftrightarrow ~~~ (Y,G) \textup{ has property } \mathcal P$.
	\end{enumerate}
	
\end{observation}
	\begin{example}\label{tehnika}
	Let $\mathcal P$ be a  property that is defined as follows: For any CR-dynamical system $(X,F)$, $$
	(X,F) \textup{ has property } \mathcal P ~~~ \Longleftrightarrow ~~~ \textup{ for all } x\in X, (x,x)\in F. 
	$$
	We show that $\mathcal P$ is a CR-dynamical property. Let $(X,F)$ and $(Y,G)$ be CR-dynamical systems that are topological conjugates and suppose that $(X,F)$ has property $\mathcal P$. Let $\varphi:X\rightarrow Y$ be a homeomorphism such that for each $(x_1,x_2)\in X\times X$, the following holds
 $$
 (x_1,x_2)\in F  \Longleftrightarrow (\varphi(x_1), \varphi(x_2))\in G.
 $$
We prove that also $(Y,G)$ has property $\mathcal P$ by showing that for each $y\in Y$, $(y,y)\in G$. Let $y\in Y$ and let $x=\varphi^{-1}(y)$. Then $(x,x)\in F$ and it follows that 
$$
(y,y)=(\varphi(x), \varphi(x))\in G.
$$ 
\end{example}
In \cite{BEK1,BEK2,BEK3,USS}, motivated by dynamical properties (such as topological entropy, minimality, transitivity and topologically mixing), several CR-dynamical properties are defined. In Sections \ref{s2} and \ref{s3}, several of such CR-dynamical properties that are motivated by the (initial) specification property are introduced. 

\begin{definition}
	Let $\mathcal R$ be a dynamical property and let $\mathcal P$ be a CR-dynamical property. We say that \emph{$\mathcal P$ generalizes  $\mathcal R$} if for each dynamical system $(X,f)$, 
		$$
	(X,f) \textup{ has property } \mathcal P ~~~ \Longleftrightarrow ~~~ (X,f) \textup{ has property } \mathcal R.
	$$
\end{definition}	
	\begin{example}
	Let $\mathcal R$ be the dynamical property that is defined as follows: For any dynamical system $(X,f)$, 
	$$
	(X,f) \textup{ has property } \mathcal R ~~~ \Longleftrightarrow ~~~ \textup{ for all } x\in X, f(x)=x. 
	$$
 Also, let $P$ be a CR-dynamical property that is defined in Example \ref{tehnika}. 	Note that $\mathcal P$ generalizes $\mathcal R$.
	\end{example}
	\begin{example}
		In \cite{BEK1,BEK4}, topological entropy for closed relations on compact metric spaces is defined. It is also proved in \cite[Theorem 3.19]{BEK1} that topological entropy for closed relations generalizes topological entropy.
	\end{example}
		\begin{example}
		In \cite{BEK2}, several variants of topological minimality for closed relations on compact metric spaces are defined. It follows from their definitions that each of the variants generalizes topological minimality.
	\end{example}
	
	\begin{example}
		In \cite{BEK3,USS}, several variants of topological transitivity for closed relations on compact metric spaces are defined. It follows from their definitions that each of the variants generalizes topological transitivity.
	\end{example}

	In \cite{BEK1}, a CR-dynamical property generalizing topological entropy is introduced. In \cite{BEK2}, several CR-dynamical properties extending the notion of topological minimality are developed and analyzed, while \cite{BEK3} focuses on CR-dynamical generalizations of topological transitivity. In a similar spirit, Sections \ref{s2} and \ref{s3} introduce several CR-dynamical properties that generalize the (initial) specification property, and their relationships are studied in detail in Section \ref{s4}.

When defining a closed relation $F$ on a compact metric space $X$, we not only obtain a CR-dynamical system $(X, F)$, but also two associated topological dynamical systems: the forward system $(X_F^+, \sigma_F^+)$ and the two-sided system $(X_F, \sigma_F)$. These systems are intimately connected to the original CR-dynamical system and provide additional insight into its behavior. For precise definitions, see Definitions \ref{luhca1}, \ref{shit}, and \ref{luhca2}.
	
\begin{definition}\label{luhca1}
Let $(X,F)$ be a CR-dynamical system. 
We call
$$
X_F^+=\Big\{(x_1,x_2,x_3,\ldots )\in \prod_{{ i={1}}}^{\infty}X \ | \ \textup{ for each positive integer } i, (x_{i},x_{i+1})\in {F}\Big\}
$$
\emph{ the  Mahavier product of ${F}$}, and we call
$$
X_F=\Big\{(\ldots,x_{-3},x_{-2},x_{-1},{x_0}{ ;}x_1,x_2,x_3,\ldots )\in \prod_{i={-\infty}}^{\infty}X \ | \ \textup{ for each  integer } i, (x_{i},x_{i+1})\in {F}\Big\}
$$
\emph{the two-sided  Mahavier product of ${F}$}.
\end{definition}

\begin{definition}\label{shit}
Let $(X,F)$ be a CR-dynamical system. The function  $\sigma_F^{+} : {X_F^+} \rightarrow {X_F^+}$, 
 defined by 
$$
\sigma_F^{+} ({x_1,x_2,x_3,x_4},\ldots)=({x_2,x_3,x_4},\ldots)
$$
for each $({x_1,x_2,x_3,x_4},\ldots)\in {X_F^+}$, 
is called \emph{   the shift map on ${X_F^+}$}. The function  $\sigma_F : {X_F} \rightarrow {X_F}$, 
 defined by 
$$
\sigma_F (\ldots,x_{-3},x_{-2},x_{-1},{x_0};x_1,x_2,x_3,\ldots )=(\ldots,x_{-2},x_{-1},x_{0},{x_1};x_2,x_3,x_4,\ldots )
$$
for each $(\ldots,x_{-3},x_{-2},x_{-1},{x_0};x_1,x_2,x_3,\ldots )\in {X_F}$, 
is called \emph{   the shift map on ${X_F}$}.    
\end{definition}
We use $p_1$ and $p_2$ to denote the standard projections from $X\times X$ to $X$:
$$
p_1(x,y)=x
$$
and
$$
p_2(x,y)=y
$$
for each $(x,y)\in X\times X$. 
\begin{observation}\label{juju}
Let $(X,F)$ be a CR-dynamical system such that $p_1(F)=p_2(F)=X$. Note that 
\begin{enumerate}
	\item $\sigma_F$ is a homeomorphism, 
	\item $\sigma_F^+$ is a continuous surjection, and
	\item $\sigma_F^+$ is a homeomorphism if and only if $F^{-1}$ is a restriction of a continuous function on $X$ (here, $F^{-1}=\{(y,x)\in X\times X \ | \ (x,y)\in F\}$).
\end{enumerate}

\end{observation}
\begin{definition}\label{luhca2}
	Let $(X,F)$ be a CR-dynamical system.  The dynamical system 
	\begin{enumerate}
		\item $(X_F^{+},\sigma_F^+)$ is called \emph{a Mahavier dynamical system}.
		\item $(X_F,\sigma_F)$ is called \emph{a two-sided Mahavier dynamical system}.
	\end{enumerate}
\end{definition}

Theorem \ref{povezava} gives a connection between two-sided Mahavier products $X_F$ and inverse limits $\varprojlim(X_F^{+},\sigma_F^+)$.

\begin{theorem}\label{povezava}
Let $X$ be a compact metric space and let $F$ be a closed relation on $X$. Then the following hold.
\begin{enumerate}
	\item $\varprojlim(X_F^{+},\sigma_F^+)$ is homeomorphic to $X_F$. 
	\item $(X_F,\sigma_F^{-1})$ and $(\varprojlim(X_F^{+},\sigma_F^+),\sigma)$ are topological conjugates.
\end{enumerate}  
\end{theorem}
\begin{proof}
See \cite[Theorem 4.1]{BE}.	
\end{proof}
We continue by stating and proving the following theorem.
\begin{theorem}
	Let $\left(X,F\right)$ be a CR-dynamical system. If $p_1(F)=p_2(F)=X$, then the following statements are equivalent.
	\begin{enumerate}
		\item\label{au} The dynamical system $\left(X_{F}^+,\sigma_F^+\right)$ has the specification property. 
		\item\label{bu} The dynamical system $\left(X_{F},\sigma_F\right)$ has the specification property.
		\item\label{au1} The dynamical system $\left(X_{F}^+,\sigma_F^+\right)$ has the initial specification property. 
		\item\label{bu1} The dynamical system $\left(X_{F},\sigma_F\right)$ has the initial specification property.
	\end{enumerate} 
\end{theorem}	
\begin{proof}
First, we prove the implication from \ref{au} to \ref{bu}. Suppose that $\left(X_{F}^+,\sigma_F^+\right)$ has the specification property. Let $\sigma:\varprojlim(X_F^+,\sigma_F^+)\rightarrow \varprojlim(X_F^+,\sigma_F^+)$ be the shift map on $\varprojlim(X_F^+,\sigma_F^+)$. By Theorem  \ref{ivan}, the dynamical system  $(\varprojlim(X_F^+,\sigma_F^+),\sigma)$ has the specification property. By Theorem \ref{povezava}, $(X_F,\sigma)$ has the specification property. It follows from Theorem \ref{goran} that $(X_F,\sigma_F^{-1})$ has the specification property. This proves the implication from \ref{au} to \ref{bu}. Next, we prove the implication from \ref{bu} to \ref{au}. Suppose that $\left(X_{F},\sigma_F\right)$ has the specification property. By Theorem \ref{goran}, $\left(X_{F},\sigma_F^{-1}\right)$ has the specification property and by Theorem \ref{povezava},  $(\varprojlim(X_F^+,\sigma_F^+),\sigma)$ has the specification property. It follows from Theorem \ref{ivan} that $\left(X_{F}^+,\sigma_F^+\right)$ has the specification property.

Note that we proved that \ref{au} is equivalent to \ref{bu}. Since $\sigma_F^+$ is surjective, it follows from Theorem \ref{lala} that \ref{au1} is equivalent to \ref{au}. Note that $\sigma_F$ is also surjective. Therefore, it follows from Theorem \ref{lala} that \ref{bu} is equivalent to \ref{bu1}.
\end{proof}	
\begin{observation}
		Let $\left(X,F\right)$ be a CR-dynamical system such that $X_F\neq \emptyset$. Note that $\sigma_F$ is a surjection, even in the case where $p_1(F)=p_2(F)=X$ is not true. Therefore, the following statements are equivalent.
	\begin{enumerate}
		\item\label{bucc} The dynamical system $\left(X_{F},\sigma_F\right)$ has the specification property.
		\item\label{bu1cc} The dynamical system $\left(X_{F},\sigma_F\right)$ has the initial specification property.
	\end{enumerate}
	\end{observation}
	Definition \ref{Ivansito} introduces the notions \emph{powered by} and \emph{fully powered by} a CR-dynamical property $P$ to formally capture the idea that a classical dynamical system can inherit its dynamical behavior from a closed relation via the Mahavier constructions. This provides a natural way to relate properties of CR-dynamical systems to properties of their associated (forward or two-sided) Mahavier dynamical systems.
\begin{definition}\label{Ivansito}
	Let $\mathcal P$ be a CR-dynamical property and let $(Y,g)$ be a dynamical system. We say that
	\begin{enumerate}
		\item $(Y,g)$ \emph{is powered by} $\mathcal P$, if there is a CR-dynamical system $(X,F)$ such that
		\begin{enumerate}
		\item $(X,F)$ has property $\mathcal P$, and
		\item $(Y,g)$ and $(X_F^+,\sigma_F^+)$ are topological conjugates.
		\end{enumerate}
		\item $(Y,g)$ is \emph{fully powered by} $\mathcal P$, if there is a CR-dynamical system $(X,F)$ such that
		\begin{enumerate}
		\item $(X,F)$ has property $\mathcal P$, and
		\item $(Y,g)$ and $(X_F,\sigma_F)$ are topological conjugates.
		\end{enumerate}
	\end{enumerate} 
\end{definition}
\begin{example}
	Let $A$ be an arc and let $f:A\rightarrow A$ be defined by $f(x)=x$ for any $x\in A$. We show that $(A,f)$ is powered by $\mathcal P$, where $\mathcal P$ is the CR-dynamical property that is defined in Example \ref{tehnika}. Let $(X,F)$ be a CR-dynamical system that is defined as follows: $(X,F)=([0,1],g)$, where $g:[0,1]\rightarrow [0,1]$ is the identity function. Note that $X_F^+$ is an arc and let $\varphi:A\rightarrow X_F^+$ be a homeomorphism. Note that for each $(x,y)\in A\times A$, the following holds
 $$
 (x,y)\in f  \Longleftrightarrow (\varphi(x), \varphi(y))\in F.
 $$
It follows that   
		\begin{enumerate}
		\item $(X,F)$ has property $\mathcal P$, and
		\item $(A,f)$ and $(X_F^+,\sigma_F^+)$ are topological conjugates.
		\end{enumerate}
Therefore, the dynamical system  $(A,f)$ is powered by $\mathcal P$. A similar argument can be used to see that $(A,f)$ is fully powered by $\mathcal P$.
	\end{example}
	
		\begin{example}
		Let $\mathcal P$ be any CR-dynamical property and let the CR-dynamical system $(X,F)$ have the property $\mathcal P$. Then $(X_F^+,\sigma_F^+)$ is powered by $\mathcal P$, and $(X_F,\sigma_F)$ is fully powered by $\mathcal P$.

		For example, let $I=[0,1]$ and let $F$ be the relation on $I$ as defined in \cite[Example 4.12]{BEK4}, and let $\mathcal P$ be the property that is defined by:
		For each CR-dynamical system $(X,F)$, $(X,F)$ has property $\mathcal P$ if and only if the entropy of $F$ is not equal to zero. Then $(I_F^+,\sigma_F^+)$ is powered by $\mathcal P$, and $(I_F,\sigma_F)$ is fully powered by $\mathcal P$.
			\end{example}

In Section \ref{s4}, several dynamical systems that are powered by a CR-dynamical property are presented. 
\begin{theorem}
	Let $\mathcal P$ be a CR-dynamical property and let $\mathcal R$ be the property that is defined as follows: For any dynamical system $(X,f)$, 
	$$
	(X,f) \textup{ has property } \mathcal R ~~~ \Longleftrightarrow ~~~ (X,f) \textup{ is (fully) powered by } \mathcal P. 
	$$
Then $\mathcal R$ is a dynamical property.
\end{theorem}
\begin{proof}
	Let $(X,f)$ and $(Y,g)$ be dynamical systems that are topological conjugates and suppose that $(X,f)$ has property $\mathcal R$. Also, let $\varphi:X\rightarrow Y$ be a homeomorphism such that $\varphi \circ f=g\circ \varphi$.  First, assume that $(X,f)$ is powered by  $\mathcal P$ and let $(Z,F)$ be a CR-dynamical system such that
		\begin{enumerate}
		\item $(Z,F)$ has property $\mathcal P$, and
		\item $(X,f)$ and $(Z_F^+,\sigma_F^+)$ are topological conjugates.
		\end{enumerate}
Let $\psi : Z_F^+\rightarrow X$ be a homeomorphism such that $f\circ \psi =\psi\circ \sigma_F^+$. Then $\varphi\circ \psi :Z_F^+\rightarrow Y$ is a homeomorphism such that 
$$
g\circ (\varphi\circ \psi)=(\varphi\circ \psi)\circ \sigma_F^+.
$$
It follows that 
\begin{enumerate}
		\item $(Z,F)$ has property $\mathcal P$, and
		\item $(Y,g)$ and $(Z_F^+,\sigma_F^+)$ are topological conjugates.
		\end{enumerate}
		Therefore, $(Y,g)$ is powered by  $\mathcal P$ and, hence, $(Y,g)$ has property $\mathcal R$. The second part of the proof, where we assume that $(X,f)$ is fully powered by  $\mathcal P$ is analogous to the first part of the proof. We leave the details to the reader.
\end{proof}

\begin{definition}
	Let $\mathcal P$ be a CR-dynamical property and let $\mathcal R$ be the property that is defined as follows:
	\begin{enumerate}
		\item For any dynamical system $(X,f)$, 
	$$
	(X,f) \textup{ has property } \mathcal R ~~~ \Longleftrightarrow ~~~ (X,f) \textup{ is powered by } \mathcal P. 
	$$
Then we write $\mathcal R=\textup{Power}(\mathcal P)$.
\item For any dynamical system $(X,f)$, 
	$$
	(X,f) \textup{ has property } \mathcal R ~~~ \Longleftrightarrow ~~~ (X,f) \textup{ is fully powered by } \mathcal P. 
	$$
Then we write $\mathcal R=\textup{FPower}(\mathcal P)$.
	\end{enumerate} 
\end{definition}

\begin{observation}
	Let $\mathcal P$ be a CR-dynamical property and let $(X,F)$ be a CR-dynamical system that has property $\mathcal P$. Note that it follows that
	\begin{enumerate}
		\item $(X_F^+,\sigma_F^+)$ is powered by the CR-property $\mathcal P$, i.e., $(X_F^+,\sigma_F^+)$ has the dynamical property $\textup{Power}(\mathcal P)$.
		\item $(X_F,\sigma_F)$ is fully powered by the CR-property $\mathcal P$, i.e., $(X_F,\sigma_F)$ has the dynamical property $\textup{FPower}(\mathcal P)$.
	\end{enumerate}
\end{observation}

\section{Specification-type properties for Mahavier dynamical systems}\label{s2}
In this section, we generalize the specification property from topological dynamical systems to CR-dynamical systems. Then we study dynamical systems that are (fully) powered by these properties. 

In Sections \ref{s2}, \ref{s3} and \ref{s4}, the Hausdorff metric is used, therefore, we begin this section by defining it.
\begin{definition}
Let $(X,d)$ be a compact metric space. Then we define \emph{$2^X$} by 
$$
2^{X}=\{A\subseteq X \ | \ A \textup{ is a non-empty closed subset of } X\}.
$$
Let $\varepsilon >0$ and let $A\in 2^X$. Then we define  \emph{$N_d(\varepsilon,A)$} by 
$$
N_d(\varepsilon,A)=\bigcup_{a\in A}B(a,\varepsilon).
$$
The function \emph{$H_d:2^X\times 2^X\rightarrow \mathbb R$}, defined by
$$
H_d(A,B)=\inf\{\varepsilon>0 \ | \ A\subseteq N_d(\varepsilon,B), B\subseteq N_d(\varepsilon,A)\},
$$
for all  $A,B\in 2^X$, is called \emph{the Hausdorff metric} on $2^X$.  
\end{definition}
Let $(X,d)$ be a compact metric space. The Hausdorff metric on $2^X$ is in fact a metric on $2^X$ and the metric space $(2^X,H_d)$ is called \emph{the hyperspace of the space $(X,d)$}. For more informaton on the topic, see \cite{nadler}.  We continue with CR-dynamical systems.

\begin{definition}
Let $(X,F)$ be a CR-dynamical system, let $x\in X$, and let $k,\ell$ be non-negative integers such that $k\leq \ell$. 
 We use $F^{[k,\ell]}(x)$ to denote  
 $$
 F^{[k,\ell]}(x)=\Big(F^k(x),F^{k+1}(x),F^{k+2}(x),\ldots,F^{\ell}(x)\Big).
 $$
  We say that $F^{[k,\ell]}(x)$ is \emph{the $[k,\ell]$-orbit segment of the point $x$}. 

\end{definition}

\begin{definition}
Let $(X,F)$ be a CR-dynamical system, let $n$ be a positive integer and for each $j\in \{1,2,3,\ldots,n\}$, let 
\begin{enumerate}
	\item $k_j$ and $\ell_j$ be non-negative integers such that $k_j\leq \ell_j$, and
	\item $x_j\in X$.
\end{enumerate}
We say that the $n$-tuple 
$$
\Bigg(F^{[k_1,\ell_1]}(x_1),F^{[k_2,\ell_2]}(x_2),F^{[k_3,\ell_3]}(x_3),\ldots ,F^{[k_n,\ell_n]}(x_n)\Bigg)
$$
 is \emph{an $n$-specification} or just \emph{a specification in $(X,F)$}.
\end{definition}
\begin{definition}
Let $(X,F)$ be a CR-dynamical system, let $n$ and $N$ be positive integers, and for each $j\in \{1,2,3,\ldots,n\}$, let 
\begin{enumerate}
	\item $k_j$ and $\ell_j$ be non-negative integers such that $k_j\leq \ell_j$, and
	\item $x_j\in X$.
\end{enumerate}
and let 
 $$
 \mathcal S=\Bigg(F^{[k_1,\ell_1]}(x_1),F^{[k_2,\ell_2]}(x_2),F^{[k_3,\ell_3]}(x_3),\ldots ,F^{[k_n,\ell_n]}(x_n)\Bigg)
 $$
  be a specification in $(X,F)$. Then we say that the specification $\mathcal S$ is \emph{an $N$-spaced specification}, if for each $j\in\{1,2,3,\ldots,n-1\}$,
$$
k_{j+1}-\ell_j\geq N.
$$
\end{definition}
\begin{definition}
Let $(X,F)$ be a CR-dynamical system, let $d$ be the metric on $X$, let $H_d$ be the Hausdorff metric on $2^X$, let $N$ be a positive integer, let $\varepsilon>0$, let $y\in X$, and let 
$$
\mathcal S=\Bigg(F^{[k_1,\ell_1]}(x_1),F^{[k_2,\ell_2]}(x_2),F^{[k_3,\ell_3]}(x_3),\ldots ,F^{[k_n,\ell_n]}(x_n)\Bigg)
$$
be an $N$-spaced specification in $(X,F)$. We say that
\begin{enumerate}
	\item the specification $\mathcal S$ is \emph{$\varepsilon$-traced in $(X,F)$ by $y$} if for each $i\in \{1,2,3,\ldots,n\}$ and for each $j\in \{k_i,k_i+1,k_i+2,\ldots,\ell_i\}$,
$$
d(F^j(y),F^j(x_i))\leq \varepsilon.
$$
\item the specification $\mathcal S$ is \emph{Hausdorff  $\varepsilon$-traced in $(X,F)$ by $y$} if for each $i\in \{1,2,3,\ldots,n\}$ and for each $j\in \{k_i,k_i+1,k_i+2,\ldots,\ell_i\}$,
$$
H_d(F^j(y),F^j(x_i))\leq \varepsilon.
$$
\end{enumerate} 
\end{definition}
\begin{definition}
Let $(X,F)$ be a CR-dynamical system. We say that $(X,F)$ has
\begin{enumerate}
	\item \emph{the specification property} (or $\mathcal{SP}$) if for each $\varepsilon >0$, there is a positive integer $N$ such that for any $N$-spaced specification $\mathcal S$ in $(X,F)$, there is $y\in X$ such that $\mathcal S$ is $\varepsilon$-traced in $(X,F)$ by $y$.
\item \emph{the Hausdorff specification property} (or $\mathcal{HSP}$) if for each $\varepsilon >0$, there is a positive integer $N$ such that for any $N$-spaced specification $\mathcal S$ in $(X,F)$, there is $y\in X$ such that $\mathcal S$ is Hausdorff  $\varepsilon$-traced in $(X,F)$ by $y$.
\end{enumerate} 
\end{definition}
\begin{observation}
	Note that if $(X,F)$ is a CR-dynamical system such that for some $x_0\in X$, 
	$$
	X\times \{x_0\}\subseteq F,
	$$
	Then for any $x,y\in X$, $d(F(x),F(y))=0$, and, therefore, $(X,F)$ has the specification property.
\end{observation}
\begin{theorem}\label{klavz}
	Let $\mathcal P\in \{\mathcal{SP},\mathcal{HSP}\}$. Then the following hold.
	\begin{enumerate}
		\item $\mathcal P$ is a CR-dynamical property that generalizes the specification property.
		\item For any dynamical system $(X,f)$, 
		$$
		(X,f) \textup{ has the specification property } ~~~ \Longrightarrow ~~~ (X,f) \textup{ has } \textup{Power}(\mathcal P).    
		$$ 
	\end{enumerate}
\end{theorem}
\begin{proof}
	 First, we prove that $\mathcal P$ is a CR-dynamical property that generalizes the specification property. To see that $\mathcal P$ is a CR-dynamical property, let $(X,F)$ and $(Y,G)$ be CR-dynamical systems that are topologically conjugate and let $\varphi:X\rightarrow Y$ be a homeomorphism such that for each $(x,y)\in X\times X$, the following holds
 $$
 (x,y)\in F  \Longleftrightarrow (\varphi(x), \varphi(y))\in G.
 $$
Let $d_X$ be the metric on $X$ and $d_Y$ the metric on $Y$ and suppose that $(X,F)$ has property $\mathcal P$. We prove  that $(Y,G)$ has property $\mathcal P$ by considering the following cases. 
	\begin{enumerate}
		\item $\mathcal P=\mathcal{SP}$.  We show that for each $\varepsilon >0$, there is a positive integer $N$ such that for any $N$-spaced specification $\mathcal S$ in $(Y,G)$, there is $y\in Y$ such that $\mathcal S$ is $\varepsilon$-traced in $(Y,G)$ by $y$. Let $\varepsilon >0$ and let $\delta>0$ be such that for all $x_1,x_2\in X$,
$$
d_X(x_1,x_2)<\delta ~~~ \Longrightarrow ~~~ d_Y(\varphi(x_1),\varphi(x_2))<\varepsilon.
$$ 
Since $(X,F)$ has property $\mathcal P$, there is a positive integer $N$ such that for any $N$-spaced specification $\mathcal S$ in $(X,F)$, there is $x\in X$ such that $\mathcal S$ is $\frac{\delta}{2}$-traced in $(X,F)$ by $x$. Choose and fix such a positive integer $N$ and let 
$$
\mathcal S=\Bigg(G^{[k_1,\ell_1]}(y_1),G^{[k_2,\ell_2]}(y_2),G^{[k_3,\ell_3]}(y_3),\ldots ,G^{[k_n,\ell_n]}(y_n)\Bigg)
$$
 be an $N$-spaced specification in $(Y,G)$.  Then 
 $$
\mathcal C=\Bigg(F^{[k_1,\ell_1]}(\varphi^{-1}(y_1)),F^{[k_2,\ell_2]}(\varphi^{-1}(y_2)),F^{[k_3,\ell_3]}(\varphi^{-1}(y_3)),\ldots ,F^{[k_n,\ell_n]}(\varphi^{-1}(y_n))\Bigg)
$$
 is an $N$-spaced specification in $(X,F)$. Let $x$ be such that $\mathcal C$ is $\frac{\delta}{2}$-traced in $(X,F)$ by $x$ and let $y=\varphi(x)$. To prove that $\mathcal S$ is $\varepsilon$-traced in $(Y,G)$ by $y$, let $i\in \{1,2,3,\ldots,n\}$ and let $j\in \{k_i,k_i+1,k_i+2,\ldots,\ell_i\}$. We prove that $d_Y(G^j(y),G^j(y_i))\leq \varepsilon$. First, note that 
 $d_X(F^j(x),F^j(\varphi^{-1}(y_i)))\leq\frac{\delta}{2}$. Let $t_1\in F^j(x)$ and let $t_2\in F^j(\varphi^{-1}(y_i))$ be such that $d_X(t_1,t_2)\leq \frac{\delta}{2}$. Then $\varphi(t_1)\in \varphi(F^j(x))$ and $\varphi(t_2)\in \varphi(F^j(\varphi^{-1}(y_i)))$.   By Theorem \ref{dekanicaMaja}, $\varphi(t_1)\in G^j(\varphi(x))$ and $\varphi(t_2)\in G^j(\varphi(\varphi^{-1}(y_i)))$. Note that $G^j(\varphi(x))=G^j(y)$, that $G^j(\varphi(\varphi^{-1}(y_i)))=G^j(y_i)$, and that  $d_Y(\varphi(t_1),\varphi(t_2))<\varepsilon$. Therefore, 
 $$
 d_Y(G^j(y),G^j(y_i))\leq \varepsilon.
 $$ 
 This proves that $\mathcal P$ is a CR-dynamical property. It follows from the definition of property $\mathcal{SP}$ that $\mathcal P$ generalizes the specification property.
 \item $\mathcal P=\mathcal{HSP}$.  We show  that for each $\varepsilon >0$, there is a positive integer $N$ such that for any $N$-spaced specification $\mathcal S$ in $(Y,G)$, there is $y\in Y$ such that $\mathcal S$ is Hausdorff $\varepsilon$-traced in $(Y,G)$ by $y$. Let $\varepsilon >0$ and let $\delta>0$ be such that for all $x_1,x_2\in X$,
$$
d_X(x_1,x_2)<\delta ~~~ \Longrightarrow ~~~ d_Y(\varphi(x_1),\varphi(x_2))<\varepsilon.
$$ 
Since $(X,F)$ has property $\mathcal P$, there is a positive integer $N$ such that for any $N$-spaced specification $\mathcal S$ in $(X,F)$, there is $x\in X$ such that $\mathcal S$ is Hausdorff $\frac{\delta}{2}$-traced in $(X,F)$ by $x$. Choose and fix such a positive integer $N$ and let 
$$
\mathcal S=\Bigg(G^{[k_1,\ell_1]}(y_1),G^{[k_2,\ell_2]}(y_2),G^{[k_3,\ell_3]}(y_3),\ldots ,G^{[k_n,\ell_n]}(y_n)\Bigg)
$$
 be an $N$-spaced specification in $(Y,G)$.  Then 
 $$
\mathcal C=\Bigg(F^{[k_1,\ell_1]}(\varphi^{-1}(y_1)),F^{[k_2,\ell_2]}(\varphi^{-1}(y_2)),F^{[k_3,\ell_3]}(\varphi^{-1}(y_3)),\ldots ,F^{[k_n,\ell_n]}(\varphi^{-1}(y_n))\Bigg)
$$
 is an $N$-spaced specification in $(X,F)$. Let $x$ be such that $\mathcal C$ is Hausdorff $\frac{\delta}{2}$-traced in $(X,F)$ by $x$ and let $y=\varphi(x)$. To prove that $\mathcal S$ is Hausdorff $\varepsilon$-traced in $(Y,G)$ by $y$, let $i\in \{1,2,3,\ldots,n\}$ and let $j\in \{k_i,k_i+1,k_i+2,\ldots,\ell_i\}$. We prove that $H_{d_Y}(G^j(y),G^j(y_i))\leq \varepsilon$ by showing that 
 \begin{enumerate}
 \item for each $z\in G^j(y)$, there is $w\in G^j(y_i)$ such that $d_Y(z,w)\leq \varepsilon$, and 
 \item for each $z\in G^j(y_i)$, there is $w\in G^j(y)$ such that $d_Y(z,w)\leq \varepsilon$. 
 \end{enumerate} 
 First, let $z\in G^j(y)$. Then  $\varphi^{-1}(z)\in \varphi^{-1}(G^j(y))$ and, and it follows from  Theorem \ref{tehnika} that $\varphi^{-1}(G^j(y))\subseteq F^j(\varphi^{-1}(y))$. Therefore, $\varphi^{-1}(z)\in F^j(x)$. Since $H_{d_X}(F^j(x),F^j(\varphi^{-1}(y_i)))\leq\frac{\delta}{2}$, it follows that there is $t\in F^j(\varphi^{-1}(y_i))$ such that $d_X(\varphi^{-1}(z),t)<\frac{\delta}{2}$. Choose and fix such a point $t$ and let $w=\varphi(t)$. Therefore,  $d_y(z,w)<\varepsilon$.  
 
 Next, let $z\in G^j(y_i)$. Then  $\varphi^{-1}(z)\in \varphi^{-1}(G^j(y_i))$ and, it follows from Theorem \ref{tehnika} that $\varphi^{-1}(G^j(y_i))\subseteq F^j(\varphi^{-1}(y_i))$. Therefore, $\varphi^{-1}(z)\in F^j(\varphi^{-1}(y_i))$. Again, it follows that there is $t\in F^j(x)$ such that $d_X(\varphi^{-1}(z),t)<\frac{\delta}{2}$. Choose and fix such a point $t$ and let $w=\varphi(t)$. Therefore,  $d_y(z,w)<\varepsilon$.  This proves that  $H_{d_Y}(G^j(y),G^j(y_i))\leq \varepsilon$ and, therefore, $\mathcal P$ is a CR-dynamical property. It follows from the definition of property $\mathcal{HSP}$ that $\mathcal P$ generalizes the specification property.
\end{enumerate} 
This completes the proof of the first part of the theorem. Next, we prove that for any dynamical system $(X,f)$, 
		$$
		(X,f) \textup{ has the specification property } ~~~ \Longrightarrow ~~~ (X,f) \textup{ has } \textup{Power}(\mathcal P).    
		$$ 
		Let $(X,f)$ be a dynamical system and suppose that $(X,f)$ has the specification property. We show that $(X,f)$ has the property $\textup{Power}(\mathcal P)$ by showing that there is a CR-dynamical system $(Z,F)$ such that 
		\begin{enumerate}
		\item $(Z,F)$ has property $\mathcal P$, and
		\item $(X,f)$ and $(Z_F^+,\sigma_F^+)$ are topological conjugates.
		\end{enumerate}
		Let $(Z,F)=(X,f)$. Note that $(Z,F)$ has property $\mathcal P$. Also, let $\varphi:X\rightarrow Z_F^+$ be defined by
		$$
		\varphi(x)=(x,f(x),f^2(x),f^3(x),\ldots)
		$$
		for each $x\in X$.  Then $\varphi\circ f=\sigma_F^+\circ \varphi$ and it follows that $(X,f)$ and $(Z_F^+,\sigma_F^+)$ are topological conjugates. 
		\end{proof}
		\begin{observation}
		Let $\mathcal P\in \{\mathcal{SP},\mathcal{HSP}\}$. Note that for any dynamical system $(X,f)$, where $f$ is a homeomorphism, 
		$$
		(X,f) \textup{ has the specification property } ~~~ \Longrightarrow ~~~ (X,f) \textup{ has } \textup{FPower}(\mathcal P).    
		$$ 
		\end{observation}
		
		\begin{theorem}\label{Hausdorff}
	Let $(X,F)$ be a CR-dynamical system. If $(X,F)$ has the Hausdorff specification property, then $(X,F)$ has the specification property.
\end{theorem}
\begin{proof}
	Suppose that $(X,F)$ has the Hausdorff specification property. To prove that $(X,F)$ has the specification property, let $\varepsilon >0$ and let $N$ be a positive integer such that for any specification in $(X,F)$, there is a point $y\in X$ such that $\mathcal S$ is Hausdorff $\varepsilon$-traced in $(X,F)$ by $y$. We prove that for any $N$-spaced  specification $\mathcal S$ in $(X,F)$, there is $y\in X$ such that $\mathcal S$ is $\varepsilon$-traced in $(X,F)$ by $y$. To do so, let 
$$
	\mathcal S=\Big((x_{k_1}^{1},x_{k_1+1}^{1},x_{k_1+2}^{1},\ldots,x_{\ell_1}^{1}),\ldots ,(x_{k_n}^{n},x_{k_n+1}^{n},x_{k_n+2}^{n},\ldots,x_{\ell_n}^{n})\Big)
	$$
	 be an $N$-spaced specification in $(X,F)$ and let $y\in X$ be such that $\mathcal S$ is Hausdorff $\varepsilon$-traced in $(X,F)$ by $y$. To see that $\mathcal S$ is $\varepsilon$-traced in $(X,F)$ by $y$, let $i\in \{1,2,3,\ldots,n\}$ and let $j\in \{k_i,k_i+1,k_i+2,\ldots,\ell_i\}$. Then 
$$
d(F^j(y),F^j(x_i))\leq H_d(F^j(y),F^j(x_i))\leq \varepsilon.
$$
\end{proof}

		\begin{corollary}
		Let $(X,f)$ be a dynamical system.	Note that
		\begin{enumerate}
			\item if $(X,f)$ has $\textup{Power}(\mathcal{HSP})$, then $(X,f)$ has $\textup{Power}(\mathcal{SP})$.
			\item if $(X,f)$ has $\textup{FPower}(\mathcal{HSP})$, then $(X,f)$ has $\textup{FPower}(\mathcal{SP})$.
		\end{enumerate} 
		\end{corollary}
		\begin{proof}
			The corollary follows directly from Theorem \ref{Hausdorff}.
		\end{proof}
		
In the following example, we give a CR-dynamical system $(X,F)$ that has the specification property but does not have the Hausdorff specification property.	
	\begin{example}\label{monica}
		Let $X = [0,1]$ and let
$$
F = \left( \left[0, \frac{1}{2} \right] \times \{0\} \right) \cup \left( \left[ \frac{1}{2}, 1 \right] \times \{1\} \right) \cup (\{1\} \times [0,1]).
$$   
To prove that $(X,F)$ has   the specification property, let $\varepsilon>0$. Next, let $N=5$, let  
	$$
	\mathcal S=\Bigg(F^{[k_1,\ell_1]}(x_1),F^{[k_2,\ell_2]}(x_2),F^{[k_3,\ell_3]}(x_3),\ldots ,F^{[k_n,\ell_n]}(x_n)\Bigg)
	$$
	 be a specification in $(X,F)$, and let $y=x_1$. To see that $\mathcal S$ is $\varepsilon$-traced in $(X,F)$ by $y$, let $i\in \{1,2,3,\ldots,n\}$ and let $j\in \{k_i,k_i+1,k_i+2,\ldots,\ell_i\}$. We consider the following possible cases. 
	 \begin{enumerate}
	 	\item $j<4$. Then $j\in \{k_1,k_1+1,k_1+2,\ldots,\ell_1\}$ and it follows that $i=1$. Therefore, 
$$
d(F^j(y),F^j(x_i))=d(F^j(y),F^j(x_1))=d(F^j(x_1),F^j(x_1))=0\leq \varepsilon.
$$
\item $j\geq 4$. Then $F^j(y),F^j(x_i)\in \{\{0\},[0,1]\}$ and, therefore, 
$$
d(F^j(y),F^j(x_i))=0\leq \varepsilon.
$$
	 \end{enumerate}

Next we show that $(X,F)$ does not have the Hausdorff specification property.
Let $\varepsilon=\frac{1}{4}$ and let $N$ be any positive integer. Take the 2-specification 
$$
\mathcal S= \left( F^{[k_1,\ell_1]} (0), F^{[k_2,\ell_2]}(1) \right)
$$
 where $k_1>1$ and $k_2 - \ell_1 \geq N$. Note that $F^k(0)=\{ 0\}$ for all $k\in \mathbb N$ and $F^{\ell}(1)=[0,1]$ for all $\ell \in \mathbb N$. Let $y\in X$ be any point. We show that $\mathcal S$ is not $\varepsilon$-traced in $(X,F)$ by $y$.
We consider the following two possible cases:
\begin{enumerate}
    \item $y\in[0,\frac{1}{2})$. For all $j \in \{ k_2,k_2+1,\ldots,\ell_2\}$ we have 
    $$
    H_d(F^j(y),F^j(1))=H_d(\{0\},[0,1])=1\nleq \varepsilon.
    $$
    \item $y\in[\frac{1}{2},1].$ For all $j \in \{ k_1,k_1+1,\ldots,\ell_1\}$ we have
    $$
    H_d(F^j(y),F^j(0))=H_d([0,1], \{0\})=1\nleq \varepsilon.
    $$
\end{enumerate}
	\end{example}

\section{Initial-specification-type properties for Mahavier dynamical systems}\label{s3}
In this section, we generalize the initial specification property from topological dynamical systems to CR-dynamical systems. Then we study dynamical systems that are (fully) powered by these properties.

\begin{definition}
	Let $(X,F)$ be a CR-dynamical system, let $x\in X$, and let $\ell$ be a non-negative integer.  We say that $F^{[0,\ell]}(x)$ is \emph{an initial $\ell$-orbit segment of the point $x$}. 
\end{definition}
\begin{definition}
	Let $(X,F)$ be a CR-dynamical system, let $n$ be a positive integer, and for each $j\in \{1,2,3,\ldots,n\}$, let
	\begin{enumerate}
		\item $\ell_j$ be non-negative integers,
		\item $x_j\in X$. 
	\end{enumerate}
	We say that the $n$-tuple
	$$
	\Bigg(F^{[0,\ell_1]}(x_1),F^{[0,\ell_2]}(x_2),F^{[0,\ell_3]}(x_3),\ldots,F^{[0,\ell_n]}(x_n)\Bigg)
	$$
	 is \emph{an initial $n$-specification} or just \emph{an initial specification in $(X,F)$}.
	\end{definition}

\begin{definition}
	Let $(X,F)$ be a CR-dynamical system, let $d$ be the metric on $X$, let $H_d$ be the Hausdorff metric on $2^X$, let $n$ be a positive integer, let $m_1,m_2,m_3\ldots,m_{n-1}$ be positive integers, let $\varepsilon>0$, let $y\in X$, and let 
	$$
	\mathcal S=\Bigg(F^{[0,\ell_1]}(x_1),F^{[0,\ell_2]}(x_2),F^{[0,\ell_3]}(x_3),\ldots ,F^{[0,\ell_n]}(x_n)\Bigg)
	$$
	be an initial specification in $(X,F)$. We say that
	\begin{enumerate}
		\item the initial specification $\mathcal S$ is \emph{$\left(\varepsilon,m_1,m_2,m_3,\ldots,m_{n-1}\right)$-traced in $(X,F)$ by $y$} if for each $i\in \{1,2,3,\ldots,n\}$ and for each $j\in \{0,1,2,\ldots,\ell_i\}$,
		$$
		d(F^{\ell_1+m_1+\ell_2+m_2+\ell_3+m_3+\cdots+\ell_{i-1}+m_{i-1}+j}\left(y\right),F^j\left(x_i\right))\leq \varepsilon.
		$$
		\item the initial specification $\mathcal S$ \emph{is Hausdorff  $\left(\varepsilon,m_1,m_2,m_3,\ldots,m_{n-1}\right)$-traced in $(X,F)$ by $y$} if for each $i\in \{1,2,3,\ldots,n\}$ and for each $j\in \{0,1,2,\ldots,\ell_i\}$,
		$$
		H_d(F^{\ell_1+m_1+\ell_2+m_2+\ell_3+m_3+\cdots+\ell_{i-1}+m_{i-1}+j}\left(y\right),F^j\left(x_i\right))\leq \varepsilon.
		$$
	\end{enumerate} 
\end{definition}
\begin{definition}
	Let $(X,F)$ be a CR-dynamical system. We say that $(X,F)$ has
	\begin{enumerate}
		\item  \emph{the initial specification property} (or $\mathcal{ISP}$) if for each $\varepsilon >0$ there is a positive integer $N$ such that
		\begin{enumerate}
			\item for each positive integer $n$
			\item for all positive integers $m_1,m_2,m_3,\ldots,m_{n-1}$ such that for each $i\in\left\lbrace 1,2,\ldots,n-1\right\rbrace$, $m_i\geq N$, 		
\end{enumerate} 
		and for any initial specification $\mathcal S$ in $(X,F)$ there is $y\in X$ such that $\mathcal S$ is $\left(\varepsilon,m_1,m_2,m_3,\ldots,m_{n-1}\right)$-traced in $(X,F)$ by $y$.
		\item  \emph{the Hausdorff initial specification property} (or $\mathcal{HISP}$) if for each $\varepsilon >0$ there is a positive integer $N$ such that
		 \begin{enumerate}
			\item for each positive integer $n$
			\item for all positive integers $m_1,m_2,m_3,\ldots,m_{n-1}$ such that for each $i\in\left\lbrace 1,2,\ldots,n-1\right\rbrace$, $m_i\geq N$, 		
\end{enumerate} 
		and for any initial specification $\mathcal S$ in $(X,F)$ there is $y\in X$ such that $\mathcal S$ is Hausdorff  $\left(\varepsilon,m_1,m_2,m_3,\dots,m_{n-1}\right)$-traced in $(X,F)$ by $y$.
	\end{enumerate}
\end{definition}
\begin{theorem}\label{klavz1}
	Let $\mathcal P\in \{\mathcal{ISP},\mathcal{HISP}\}$. Then the following holds.
	\begin{enumerate}
		\item $\mathcal P$ is a CR-dynamical property that generalizes the initial specification property.
		\item For any dynamical system $(X,f)$, 
		$$
		(X,f) \textup{ has the initial specification property } ~~~ \Longrightarrow ~~~ (X,f) \textup{ has } \textup{Power}(\mathcal P).    
		$$ 
	\end{enumerate}
\end{theorem}
\begin{proof}
The proof of this theorem is analogous to the proof of Theorem \ref{klavz}. We leave the details to the reader.
\end{proof}
\begin{observation}
		Let $\mathcal P\in \{\mathcal{ISP},\mathcal{HISP}\}$. Note that for any dynamical system $(X,f)$, where $f$ is a homeomorphism, 
		$$
		(X,f) \textup{ has the initial specification property } ~~~ \Longrightarrow ~~~ (X,f) \textup{ has } \textup{FPower}(\mathcal P).    
		$$ 
		\end{observation}
		
		\begin{theorem}\label{Hausdorfff}
	Let $(X,F)$ be a CR-dynamical system. If $(X,F)$ has the Hausdorff initial specification property, then $(X,F)$ has the initial specification property.
\end{theorem}
\begin{proof}
	Suppose that $(X,F)$ has the Hausdorff initial specification property. To prove that $(X,F)$ has the initial specification property, let $\varepsilon>0$ and let $N$ be a positive integer such that for each positive integer $n$, for all positive integers $m_1,m_2,m_3,\ldots,m_{n-1}$ such that for each $i\in \{1,2,3,\ldots,n\}$, $m_i\geq N$, and for each initial $n$-specification $\mathcal S$ in $(X,F)$, there is $y\in X$ such that $\mathcal S$ is Hausdorff $(\varepsilon,m_1,m_2,m_3,\ldots,m_{n-1})$-traced in $(X,F)$ by $y$. We prove that for each positive integer $n$, for all positive integers $m_1,m_2,m_3,\ldots,m_{n-1}$ such that for each $i\in \{1,2,3,\ldots,n\}$, $m_i\geq N$, and for each initial $n$-specification $\mathcal S$, there is $y\in X$ such that $\mathcal S$ is $(\varepsilon,m_1,m_2,m_3,\ldots,m_{n-1})$-traced in $(X,F)$ by $y$.  So, let $n$ be a positive integer, let $m_1,m_2,m_3,\ldots,m_{n-1}$ be positive integers such that for each $i\in \{1,2,3,\ldots,n\}$, $m_i\geq N$, and let 
$$
	\mathcal S=\Big((x_{0}^{1},x_{1}^{1},x_{2}^{1},\ldots,x_{\ell_1}^{1}),\ldots ,(x_{0}^{n},x_{1}^{n},x_{2}^{n},\ldots,x_{\ell_n}^{n})\Big)
	$$
	 be an initial specification in $(X,F)$. Also, let $y\in X$ be such that $\mathcal S$ is Hausdorff $(\varepsilon,m_1,m_2,m_3,\ldots,m_{n-1})$-traced in $(X,F)$ by $y$. We show that the initial specification $\mathcal S$ is $(\varepsilon,m_1,m_2,m_3,\ldots,m_{n-1})$-traced in $(X,F)$ by $y$. To do so, let $i\in \{1,2,3,\ldots,n\}$ and let $j\in \{0,1,2,\ldots,\ell_i\}$. Then 
\begin{align*}
	&d(F^{\ell_1+m_1+\ell_2+m_2+\ell_3+m_3+\cdots+\ell_{i-1}+m_{i-1}+j}\left(y\right),F^j\left(x_i\right))\leq \\
	&H_d(F^{\ell_1+m_1+\ell_2+m_2+\ell_3+m_3+\cdots+\ell_{i-1}+m_{i-1}+j}\left(y\right),F^j\left(x_i\right))\leq \varepsilon.
\end{align*}
This completes the proof. 
\end{proof}

		\begin{corollary}
		Let $(X,f)$ be a dynamical system.	Note that
		\begin{enumerate}
			\item if $(X,f)$ has $\textup{Power}(\mathcal{HISP})$, then $(X,f)$ has $\textup{Power}(\mathcal{ISP})$.
			\item if $(X,f)$ has $\textup{FPower}(\mathcal{HISP})$, then $(X,f)$ has $\textup{FPower}(\mathcal{ISP})$.
		\end{enumerate} 
		\end{corollary}
		\begin{proof}
			The corollary follows directly from Theorem \ref{Hausdorfff}.
		\end{proof}
		
\section{Specification properties Versus initial specification properties}\label{s4}
In this section, we study relationships between specification properties ($\mathcal{SP}$ and $\mathcal{HSP}$) and corresponding initial specification properties ($\mathcal{ISP}$ and $\mathcal{HISP}$).
\subsection{$\mathcal{SP}$ Versus $\mathcal{ISP}$}
The following example shows that, in general, the   specification property is not equivalent to the   initial specification property. 
\begin{example}
	Let $X=[0,1]$ and let $F=[0,1]\times \{1\}$. To prove that $(X,F)$ has   the specification property, let $\varepsilon>0$. Then, let $N=1$, let  
	$$
	\mathcal S=\Bigg(F^{[k_1,\ell_1]}(x_1),F^{[k_2,\ell_2]}(x_2),F^{[k_3,\ell_3]}(x_3),\ldots ,F^{[k_n,\ell_n]}(x_n)\Bigg)
	$$
	 be an $N$-spaced specification in $(X,F)$, and let $y=x_1$. To see that $\mathcal S$ is $\varepsilon$-traced in $(X,F)$ by $y$, let $i\in \{1,2,3,\ldots,n\}$ and let $j\in \{k_i,k_i+1,k_i+2,\ldots,\ell_i\}$. We consider the following possible cases. 
	 \begin{enumerate}
	 	\item $j=0$. Then $j=k_1$ and it follows that  
$$
d(F^j(y),F^j(x_i))=d(F^{k_1}(x_1),F^{k_1}(x_1))=d(\{x_1\},\{x_1\})=0\leq \varepsilon.
$$
	 	\item $j\neq 0$. Then 
$$
d(F^j(y),F^j(x_i))=d(\{1\},\{1\})=0\leq \varepsilon.
$$
	 \end{enumerate}
Next, we show that $(X,F)$ does not have  the initial specification property. To see this, let $\varepsilon=\frac{1}{4}$ and let $N$ be any positive integer. Also, let $m_1$ be a positive integer such that $m_1\geq N$, and let $\mathcal S=\Bigg(F^{[0,1]}(1),F^{[0,1]}(0)\Bigg)$ be an initial $2$-specification in $(X,F)$. Note that in this case, $x_1=1$, $x_2=0$, and $\ell_1=\ell_2=1$. Finally, let $y\in X$ be any point. We show that $\mathcal S$ is not $(\varepsilon,m_1)$-traced in $(X,F)$ by $y$. We prove this by showing that it does not hold that for each $i\in \{1,2\}$ and for each $j\in \{0,1,2,\ldots,\ell_i\}$,
		$$
		d(F^{\ell_1+m_1+\ell_2+m_2+\ell_3+m_3+\cdots+\ell_{i-1}+m_{i-1}+j}\left(y\right),F^j\left(x_i\right))\leq \varepsilon.
		$$
Let $i=2$ and let $j=0$. Then 
\begin{align*}
		&d(F^{\ell_1+m_1+\ell_2+m_2+\ell_3+m_3+\cdots+\ell_{i-1}+m_{i-1}+j}\left(y\right),F^j\left(x_i\right))=\\
		&d(F^{\ell_1+m_1+j}\left(y\right),F^j\left(x_i\right))=d(F^{1+1+0}\left(y\right),F^0\left(0\right))=d(\{1\},\{0\})=1\not\leq \varepsilon.
\end{align*}
\end{example}

\begin{theorem}\label{prvi}
	Let $(X,F)$ be a CR-dynamical system. If there is a positive integer $n_0$ such that for all $x,y\in X$,   $F^{n_0}(x)\cap F^{n_0}(y)\neq \emptyset$, then $(X,F)$ has the specification property. 
\end{theorem}
\begin{proof}
	Let $n_0$ be a positive integer such that for all $x,y\in X$, $F^{n_0}(x)\cap F^{n_0}(y)\neq \emptyset$. Note that it follows that $p_1(X)=X$. To prove that $(X,F)$ has   the specification property, let $\varepsilon>0$. Next, let $N=n_0$, let  
	$$
	\mathcal S=\Bigg(F^{[k_1,\ell_1]}(x_1),F^{[k_2,\ell_2]}(x_2),F^{[k_3,\ell_3]}(x_3),\ldots ,F^{[k_n,\ell_n]}(x_n)\Bigg)
	$$
	 be an $N$-spaced specification in $(X,F)$, and let $y=x_1$. To see that $\mathcal S$ is $\varepsilon$-traced in $(X,F)$ by $y$, let $i\in \{1,2,3,\ldots,n\}$ and let $j\in \{k_i,k_i+1,k_i+2,\ldots,\ell_i\}$. We consider the following possible cases. 
	 \begin{enumerate}
	 	\item $j<N$. Then $j\in \{k_1,k_1+1,k_1+2,\ldots,\ell_1\}$ and it follows that $i=1$. Therefore, 
$$
d(F^j(y),F^j(x_i))=d(F^j(x_1),F^j(x_1))=0\leq \varepsilon.
$$
\item $j\geq N$. Then $F^j(y)\cap F^j(x_i)\neq \emptyset$ and, therefore, $d(F^j(y),F^j(x_i))=0\leq \varepsilon$.
	 \end{enumerate}
\end{proof}

\begin{theorem}\label{drugi}
	Let $(X,F)$ be a CR-dynamical system. If there is a positive integer $n_0$ such that for each $x\in X$,   $F^{n_0}(x)=X$, then $(X,F)$ has   the specification property as well as  the initial specification property. 
\end{theorem}
\begin{proof}
	It follows from Theorem \ref{prvi} that $(X,F)$ has   the specification property.
	 To prove that $(X,F)$ has  the initial specification property, let $\varepsilon>0$ and let $N=n_0$. Also, let $n$ be a positive integer, let $m_1,m_2,m_3,\ldots,m_{n-1}$ be positive integers such that for each $i\in \{1,2,3,\ldots ,n-1\}$, $m_i\geq N$, and let  
	$$
	\mathcal S=\Bigg(F^{[0,\ell_1]}(x_1),F^{[0,\ell_2]}(x_2),F^{[0,\ell_3]}(x_3),\ldots ,F^{[0,\ell_n]}(x_n)\Bigg)
	$$
	 be an initial specification in $(X,F)$, and let $y=x_1$. To see that the initial specification $\mathcal S$ is $(\varepsilon,m_1,m_2,m_3,\ldots,m_{n-1})$-traced in $(X,F)$ by $y$, let $i\in \{1,2,3,\ldots,n\}$ and let $j\in \{0,1,2,\ldots,\ell_i\}$. We consider the following possible cases. 
	 \begin{enumerate}
	 	\item $i=1$. Then $j\in \{0,1,2,\ldots,\ell_1\}$ and it follows that 
$$
		d(F^{\ell_1+m_1+\ell_2+m_2+\ell_3+m_3+\cdots+\ell_{i-1}+m_{i-1}+j}\left(y\right),F^j\left(x_i\right))=		d(F^{j}\left(x_1\right),F^j\left(x_1\right))\leq \varepsilon.
$$
\item $i>1$. Then $F^{\ell_1+m_1+\ell_2+m_2+\ell_3+m_3+\cdots+\ell_{i-1}+m_{i-1}+j}\left(y\right)=X$ (since $m_1\geq N=n_0$) and, therefore, 
$$
d(F^{\ell_1+m_1+\ell_2+m_2+\ell_3+m_3+\cdots+\ell_{i-1}+m_{i-1}+j}\left(y\right),F^j\left(x_i\right))=d(X,F^j\left(x_i\right))=0\leq \varepsilon
$$
since $F^j\left(x_i\right)$ is a non-empty subset of $X$. 
	 \end{enumerate}
\end{proof}
\begin{example}
	Let $X=[0,1]$ and let 
	$$
	F=\left(\left[0,\frac{1}{2}\right]\times \{0\}\right)\cup \left(\{0\}\times \left[0,\frac{1}{2}\right]\right)\cup \left(\left[\frac{1}{2},1\right]\times \{1\}\right)\cup\left(\{1\}\times \left[\frac{1}{2},1\right]\right).
	$$
	Then for each $x\in X$, $F^{4}(x)=X$. Therefore, by Theorem \ref{drugi}, $(X,F)$ has   the specification property as well as  the initial specification property. 
\end{example}
In the following example, we demonstrate that even if $p_1(F)=p_2(F)=X$, the    specification property and   the initial specification property may not be equivalent.
\begin{example}\label{ex3}
	Let $X=[0,1]$ and let 
	$$
	F=\left(\left[0,\frac{1}{2}\right]\times \{0\}\right)\cup \left(\left[\frac{1}{2},1\right]\times \{1\}\right)\cup\left(\{1\}\times \left[0,1\right]\right).
	$$
	 We already know that $(X,F)$ has   the specification property; see Example \ref{monica}. 
Next, we show that $(X,F)$ does not have  the initial specification property. To see this, let $\varepsilon=\frac{1}{8}$ and let $N$ be any positive integer. Also, let $m_1$ be a positive integer such that $m_1\geq N$, and let $\mathcal S=\Big(F^{[0,1]}(0),F^{[0,1]}(\frac{3}{4})\Big)$ be an initial $2$-specification in $(X,F)$. Note that in this case, $x_1=0$, $x_2=\frac{3}{4}$, and $\ell_1=\ell_2=1$. Finally, let $y\in X$ be any point. We show that $\mathcal S$ is not $(\varepsilon,m_1)$-traced in $(X,F)$ by $y$. We prove this by showing that it does not hold that for each $i\in \{1,2\}$ and for each $j\in \{0,1,2,\ldots,\ell_i\}$,
		$$
		d(F^{\ell_1+m_1+\ell_2+m_2+\ell_3+m_3+\cdots+\ell_{i-1}+m_{i-1}+j}\left(y\right),F^j\left(x_i\right))\leq \varepsilon.
		$$
		 We consider the following two possible cases.
\begin{enumerate}

\item $y\in [0,\frac{1}{2})$. Let $i=2$ and let $j=0$. Then 
\begin{align*}
		&d(F^{\ell_1+m_1+\ell_2+m_2+\ell_3+m_3+\cdots+\ell_{i-1}+m_{i-1}+j}\left(y\right),F^j\left(x_i\right))=\\
		&d(F^{\ell_1+m_1+j}\left(y\right),F^j\left(x_2\right))=d\left(F^{1+m_1+0}\left(y\right),F^0\left(\frac{3}{4}\right)\right)=d\left(\{0\},\left\{\frac{3}{4}\right\}\right)=\frac{3}{4}\not\leq \varepsilon.
\end{align*}
\item $y\in [\frac{1}{2},1]$. Let $i=1$ and let $j=0$. Then 
\begin{align*}
		&d(F^{\ell_1+m_1+\ell_2+m_2+\ell_3+m_3+\cdots+\ell_{i-1}+m_{i-1}+j}\left(y\right),F^j\left(x_i\right))=\\
		&d(F^{j}\left(y\right),F^j\left(x_1\right))=d\left(F^{0}\left(y\right),F^0\left(0\right)\right)=d\left(\{y\},\{0\}\right)=y\not\leq \varepsilon.
\end{align*}
\end{enumerate}
\end{example}
In Example \ref{ex3}, a CR-dynamical system $(X,F)$ is constructed in such a way that $p_1(F)=p_2(F)=X$ and such that
\begin{enumerate}
	\item $(X,F)$ has   the specification property, and
	\item $(X,F)$ does not have  the initial specification property.
\end{enumerate}
This means that   the specification property does not imply  the initial specification property (even in the case $p_1(F)=p_2(F)=X$).
In the following theorem, we prove that 
the   initial specification property implies   the specification property.
\begin{theorem}\label{tretji}
	Let $(X,F)$ be a CR-dynamical system. 
	Suppose that $(X,F)$ has the  initial specification property. Then $(X,F)$ has the specification property. 
\end{theorem}
\begin{proof} 
Let $\varepsilon > 0$. Also, let $N$ be a positive integer such that for each positive integer $n$, for all positive integers $m_1,m_2,m_3,\ldots,m_{n-1}$ such that for each $i\in \{1,2,3,\ldots,n\}$, $m_i\geq N$, and for each initial $n$-specification $\mathcal S$, there is $y\in X$ such that $\mathcal S$ is $(\varepsilon,m_1,m_2,m_3,\ldots,m_{n-1})$-traced in $(X,F)$ by $y$. Such $N$ does exist since  $(X,F)$ has  the initial specification property. We prove that for any $N$-spaced specification $\mathcal S$ in $(X,F)$, there is $y\in X$ such that $\mathcal S$ is $\varepsilon$-traced in $(X,F)$ by $y$. Let 
$$
\mathcal S=\Bigg(F^{[k_1,\ell_1]}(x_1),F^{[k_2,\ell_2]}(x_2),F^{[k_3,\ell_3]}(x_3),\ldots ,F^{[k_n,\ell_n]}(x_n)\Bigg)
$$
be an $N$-spaced specification in $(X,F)$. Also, for each $i\in \{1,2,3,\ldots, n\}$, let $z_i\in F^{k_i}(x_i)$. Then 
$$
\mathcal C=\Bigg(F^{[0,\ell_1-k_1]}(z_1),F^{[0,\ell_2-k_2]}(z_2),F^{[0,\ell_3-k_3]}(z_3),\ldots ,F^{[0,\ell_n-k_n]}(z_n)\Bigg)
$$
is an initial $n$-specification. For each $i\in \{1,2,3,\ldots, n-1\}$, let $m_i=k_{i+1}-\ell_i$. Note that for each $i\in \{1,2,3,\ldots, n-1\}$, $m_i\geq N$. Let $z\in X$ be such that $\mathcal C$ is $(\varepsilon,m_1,m_2,m_3,\ldots,m_{n-1})$-traced in $(X,F)$ by $z$. Also, let $y\in X$ be such that $z\in F^{k_1}(y)$. We claim that $\mathcal S$ is $\varepsilon$-traced in $(X,F)$ by $y$.
Let $i\in \{1,2,3,\ldots,n\}$ and let $j\in \{k_i,k_i+1,k_i+2,\ldots,\ell_i\}$. Then 
\begin{align*}
&d(F^j(y),F^j(x_i))=d(F^{j-k_1+k_1}(y),F^{j-k_i+k_i}(x_i))=d(F^{j-k_1}(F^{k_1}(y)),F^{j-k_i}(F^{k_i}(x_i)))=\\
&d(F^{(\ell_1-k_1)+m_1+(\ell_2-k_2)+m_2+(\ell_3-k_3)+m_3+\ldots+(\ell_{i-1}-k_{i-1})+m_{i-1}+(j-k_i)}(F^{k_1}(y)),F^{j-k_i}(F^{k_i}(x_i)))\leq\\
&d(F^{(\ell_1-k_1)+m_1+(\ell_2-k_2)+m_2+(\ell_3-k_3)+m_3+\ldots+(\ell_{i-1}-k_{i-1})+m_{i-1}+(j-k_i)}(z),F^{j-k_i}(z_i))\leq \varepsilon.
\end{align*}
This completes the proof. 
\end{proof}
\subsection{$\mathcal{HSP}$ Versus $\mathcal{HISP}$}

The following example shows that in general, the  Hausdorff specification property is not equivalent to the  Hausdorff initial specification property. 
\begin{example}
	Let $X=[0,1]$ and let $F=[0,1]\times \{1\}$. To prove that $(X,F)$ has the Hausdorff specification property, let $\varepsilon>0$. Then, let $N=1$, let  
	$$
	\mathcal S=\Bigg(F^{[k_1,\ell_1]}(x_1),F^{[k_2,\ell_2]}(x_2),F^{[k_3,\ell_3]}(x_3),\ldots ,F^{[k_n,\ell_n]}(x_n)\Bigg)
	$$
	 be an $N$-spaced specification in $(X,F)$, and let $y=x_1$. To see that $\mathcal S$ is Hausdorff  $\varepsilon$-traced in $(X,F)$ by $y$, let $i\in \{1,2,3,\ldots,n\}$ and let $j\in \{k_i,k_i+1,k_i+2,\ldots,\ell_i\}$. We consider the following possible cases. 
	 \begin{enumerate}
	 	\item $j=0$. Then $j=k_1$ and it follows that  
$$
H_d(F^j(y),F^j(x_i))=H_d(F^{k_1}(x_1),F^{k_1}(x_1))=H_d(\{x_1\},\{x_1\})=0\leq \varepsilon.
$$
	 	\item $j\neq 0$. Then 
$$
H_d(F^j(y),F^j(x_i))=H_d(\{1\},\{1\})=0\leq \varepsilon.
$$
	 \end{enumerate}
Next, we show that $(X,F)$ does not have the Hausdorff initial specification property. Suppose that it does have the Hausdorff initial specification property. Then $(X,F)$ has the initial specification property, which is a contradiction with Example \ref{ex3}. 
\end{example}

\begin{theorem}\label{prvi11}
	Let $(X,F)$ be a CR-dynamical system. If for each $\varepsilon>0$, there is a positive integer $n_0$ such that for all $x,y\in X$ and for each non-negative integer $j$,     $H_d(F^{n_0+j}(x),F^{n_0+j}(y))\leq \varepsilon$, then $(X,F)$ has the Hausdorff specification property. 
\end{theorem}
\begin{proof}
To prove that $(X,F)$ has the Hausdorff specification property, let $\varepsilon>0$. Let $n_0$ be a positive integer such that for all $x,y\in X$, $H_d(F^{n_0}(x),F^{n_0}(y))\leq \varepsilon$. Next, let $N=n_0$, let  
	$$
	\mathcal S=\Bigg(F^{[k_1,\ell_1]}(x_1),F^{[k_2,\ell_2]}(x_2),F^{[k_3,\ell_3]}(x_3),\ldots ,F^{[k_n,\ell_n]}(x_n)\Bigg)
	$$
	 be an $N$-spaced specification in $(X,F)$, and let $y=x_1$. To see that $\mathcal S$ is Hausdorff $\varepsilon$-traced in $(X,F)$ by $y$, let $i\in \{1,2,3,\ldots,n\}$ and let $j\in \{k_i,k_i+1,k_i+2,\ldots,\ell_i\}$. We consider the following possible cases. 
	 \begin{enumerate}
	 	\item $j<N$. Then $j\in \{k_1,k_1+1,k_1+2,\ldots,\ell_1\}$ and it follows that $i=1$. Therefore, 
$$
H_d(F^j(y),F^j(x_i))=H_d(F^j(x_1),F^j(x_1))=0\leq \varepsilon.
$$
\item $j\geq N$. Then $j\geq n_0$ and, therefore, $H_d(F^j(y),F^j(x_i))=0\leq \varepsilon$.
	 \end{enumerate}
\end{proof}	
\begin{corollary}\label{prvi1}
	Let $(X,F)$ be a CR-dynamical system. If there is a positive integer $n_0$ such that for all $x,y\in X$,   $F^{n_0}(x)=F^{n_0}(y)$, then $(X,F)$ has  the Hausdorff specification property. 
\end{corollary}
\begin{proof}
	Let $n_0$ be a positive integer such that for all $x,y\in X$, $F^{n_0}(x)=F^{n_0}(y)$. Note that it follows that for each for each $\varepsilon>0$, for all $x,y\in X$, and for each non-negative integer $j$, $H_d(F^{n_0+j}(x),F^{n_0+j}(y))=0\leq \varepsilon$.
\end{proof}
In the following example, we show that  even if there is a positive integer $n_0$ such that for each $x\in X$,   $F^{n_0}(x)=X$, $(X,F)$ does not need to have the Hausdorff initial specification property.

\begin{example}\label{exi}
	Let $X=[0,1]$ and let 
	$$
	F=\left(\left[0,\frac{1}{2}\right]\times \{0\}\right)\cup \left(\{0\}\times \left[0,\frac{1}{2}\right]\right)\cup \left(\left[\frac{1}{2},1\right]\times \{1\}\right)\cup\left(\{1\}\times \left[\frac{1}{2},1\right]\right).
	$$
	Then for each $x\in X$, $F^{4}(x)=X$. Therefore, by Corollary \ref{prvi1}, $(X,F)$ has the Hausdorff specification property. 
	
	Next, we show that $(X,F)$ does not have the Hausdorff initial specification property. To see this, let $\varepsilon=\frac{1}{4}$ and let $N$ be any positive integer. Also, let $m_1$ be a positive integer such that $m_1\geq N$, and let $\mathcal S=\Bigg(F^{[0,1]}\left(\frac{1}{4}\right),F^{[0,1]}\left(\frac{3}{4}\right)\Bigg)$ be an initial $2$-specification in $(X,F)$. Note that in this case, $x_1=\frac{1}{4}$, $x_2=\frac{3}{4}$, and $\ell_1=\ell_2=1$. Finally, let $y\in X$ be any point. We show that $\mathcal S$ is not Hausdorff $(\varepsilon,m_1)$-traced in $(X,F)$ by $y$. We prove this by showing that it does not hold that for each $i\in \{1,2\}$ and for each $j\in \{0,1,2,\ldots,\ell_i\}$,
		$$
		H_d(F^{\ell_1+m_1+\ell_2+m_2+\ell_3+m_3+\cdots+\ell_{i-1}+m_{i-1}+j}\left(y\right),F^j\left(x_i\right))\leq \varepsilon.
		$$
		We consider the following cases.
		\begin{enumerate}
			\item $y=0$. Let $i=2$ and let $j=0$. Then 
\begin{align*}
		&H_d(F^{\ell_1+m_1+\ell_2+m_2+\ell_3+m_3+\cdots+\ell_{i-1}+m_{i-1}+j}\left(y\right),F^j\left(x_i\right))=H_d(F^{\ell_1+m_1+j}\left(y\right),F^j\left(x_i\right))=\\
		&H_d(F^{1+1+0}\left(0\right),F^0\left(0\right))=H_d\left(\left[0,\frac{1}{2}\right]\cup \{1\},\{0\}\right)=1\not\leq \varepsilon.
\end{align*}
\item $y\in (0,\frac{1}{2})$. Let $i=2$ and let $j=0$. Then 
\begin{align*}
		&H_d(F^{\ell_1+m_1+\ell_2+m_2+\ell_3+m_3+\cdots+\ell_{i-1}+m_{i-1}+j}\left(y\right),F^j\left(x_i\right))=H_d(F^{\ell_1+m_1+j}\left(y\right),F^j\left(x_i\right))=\\
		&H_d(F^{1+1+0}\left(y\right),F^0\left(0\right))=H_d\left(\left[0,\frac{1}{2}\right],\{0\}\right)=\frac{1}{2}\not\leq \varepsilon.
\end{align*}
\item $y=\frac{1}{2}$. Let $i=2$ and let $j=0$. Then 
\begin{align*}
		&H_d(F^{\ell_1+m_1+\ell_2+m_2+\ell_3+m_3+\cdots+\ell_{i-1}+m_{i-1}+j}\left(y\right),F^j\left(x_i\right))=H_d(F^{\ell_1+m_1+j}\left(y\right),F^j\left(x_i\right))=\\
		&H_d\left(F^{1+1+0}\left(\frac{1}{2}\right),F^0\left(0\right)\right)=H_d\left(\left[0,1\right],\{0\}\right)=1\not\leq \varepsilon.
\end{align*}
\item $y\in (\frac{1}{2},1)$. Let $i=2$ and let $j=0$. Then 
\begin{align*}
		&H_d(F^{\ell_1+m_1+\ell_2+m_2+\ell_3+m_3+\cdots+\ell_{i-1}+m_{i-1}+j}\left(y\right),F^j\left(x_i\right))=H_d(F^{\ell_1+m_1+j}\left(y\right),F^j\left(x_i\right))=\\
		&H_d\left(F^{1+1+0}\left(y\right),F^0\left(0\right)\right)=H_d\left(\left[\frac{1}{2},1\right],\{0\}\right)=1\not\leq \varepsilon.
\end{align*}
\item $y=1$. Let $i=2$ and let $j=0$. Then 
\begin{align*}
		&H_d(F^{\ell_1+m_1+\ell_2+m_2+\ell_3+m_3+\cdots+\ell_{i-1}+m_{i-1}+j}\left(y\right),F^j\left(x_i\right))=H_d(F^{\ell_1+m_1+j}\left(y\right),F^j\left(x_i\right))=\\
		&H_d\left(F^{1+1+0}\left(1\right),F^0\left(0\right)\right)=H_d\left(\{0\}\cup \left[\frac{1}{2},1\right],\{0\}\right)=1\not\leq \varepsilon.
\end{align*}
		\end{enumerate}
\end{example}
Among other things, Example \ref{exi} demonstrates that even if $p_1(F)=p_2(F)=X$, the Hausdorff specification property and the Hausdorff initial specification property may not be equivalent. This means that the Hausdorff specification property does not imply the Hausdorff  initial specification property (even in the case $p_1(F)=p_2(F)=X$). 

We conclude the section by stating the following open problem.
\begin{problem}
	Is there a CR-dynamical system $(X,F)$ such that
	\begin{enumerate}
		\item $(X,F)$ has the Hausdorff initial specification property but it doesn't have the Hausdorff specification property. 
		\item $p_1(X)=p_2(X)=X$ and $(X,F)$ has the Hausdorff initial specification property but it doesn't have the Hausdorff specification property. 
	\end{enumerate} 
\end{problem}
\section{Acknowledgement}
This work is supported in part by the Slovenian Research Agency (research projects J1-4632, BI-HR/23-24-011, BI-US/22-24-086 and BI-US/22-24-094, and research program P1-0285). 
	

\noindent I. Bani\v c\\
              (1) Faculty of Natural Sciences and Mathematics, University of Maribor, Koro\v{s}ka 160, SI-2000 Maribor,
   Slovenia; \\(2) Institute of Mathematics, Physics and Mechanics, Jadranska 19, SI-1000 Ljubljana, 
   Slovenia; \\(3) Andrej Maru\v si\v c Institute, University of Primorska, Muzejski trg 2, SI-6000 Koper,
   Slovenia\\
             {iztok.banic@um.si}           
     
				\-
				
		\noindent G.  Erceg\\
             Faculty of Science, University of Split, Rudera Bo\v skovi\' ca 33, Split,  Croatia\\
{{gorerc@pmfst.hr}       }    

                 	\-

	\noindent I.  Jeli\' c\\
             Faculty of Science, University of Split, Rudera Bo\v skovi\' ca 33, Split,  Croatia\\
{{ivan.jelic@pmfst.hr}       }   			
		
                 	\-
					
  \noindent J.  Kennedy\\
             Department of Mathematics,  Lamar University, 200 Lucas Building, P.O. Box 10047, Beaumont, Texas 77710 USA\\
{{kennedy9905@gmail.com}       }  
\-
				

	\-
				
		%



\end{document}